\documentclass[12pt,draft]{amsart}
\usepackage{amsmath,amsfonts,amssymb,amscd,amsthm}
\newtheorem{prop}{Proposition}[section]
\newtheorem{thm}[prop]{Theorem}
\newtheorem{cor}[prop]{Corollary}
\newtheorem{ques}[prop]{Question} 
\newtheorem{adden}[prop]{Addendum} 
\newtheorem{lem}[prop]{Lemma}
\newtheorem{conj}[prop]{Conjecture}
\theoremstyle{definition}
\newtheorem{de}[prop]{Definition}
\newtheorem{example}[prop]{Example}
\newtheorem{examples}[prop]{Examples}
\theoremstyle{remark}
\newtheorem{Remarks}[prop]{Remarks}             
\def\rU{\text{\rm U}}
\def\complex{{\mathbb C}}
\def\PP{{\mathbb{P}}}
\def\que{{\mathbb Q}}
\def\real{{\mathbb R}}
\def\integer{{\mathbb Z}}
\def\zed{{\mathbb Z}}
\def\CP{{\complex\PP}}
\def\A{{\mathcal{A}}}
\def\B{{\mathcal{B}}}
\def\C{{\mathcal{C}}}
\def\J{{\mathcal{J}}}
\def\cS{{\mathcal{S}}}
\def\Aut{\mathop{\rm Aut}\nolimits}
\def\cl{\mathop{\rm cl}\nolimits}
\def\codim{\mathop{\rm codim}\nolimits}
\def\dR{\mathop{\rm dR}\nolimits}
\def\End{\mathop{\rm End}\nolimits}
\def\GL{\mathop{\rm GL}\nolimits}
\def\graph{\mathop{\rm graph}\nolimits}
\def\Hom{\mathop{\rm Hom}\nolimits}
\def\id{\mathop{\rm id}\nolimits}
\def\Int{\mathop{\rm int}\nolimits}
\def\Span{\mathop{\rm span}\nolimits}
\def\supp{\mathop{\rm supp}\nolimits}
\def\std{\mathop{\rm std}\nolimits}
\def\ostd{{\omega_{\std}}}
\begin{document}

\title{Toward a topological characterization of symplectic manifolds}
\author{Robert E. Gompf} \thanks{Partially supported by NSF grants 
DMS-9802533 and DMS-0102922.}
\address{Department of Mathematics, The University of Texas at Austin, 
Austin, TX 78712-1082}
\email{gompf@math.utexas.edu}
\date{\today}

\begin{abstract}

A fibration-like structure called a hyperpencil 
is defined on a smooth,
closed $2n$-manifold $X$, generalizing a linear system of curves on an
algebraic variety.
A deformation class of hyperpencils is shown to 
determine an isotopy class of symplectic structures on $X$. 
This provides an inverse  to Donaldson's 
program for constructing linear systems on symplectic manifolds. 
In dimensions $\le6$, work of Donaldson and Auroux provides hyperpencils
on any symplectic manifold, and the author conjectures that this extends
to arbitrary dimensions.
In dimensions where this holds, the set of deformation 
classes of hyperpencils canonically maps onto the set of isotopy classes 
of rational symplectic forms up to positive scale, topologically 
determining a dense subset of all symplectic forms up to an equivalence
relation on hyperpencils. 
In particular, the existence of a hyperpencil topologically characterizes
those manifolds in dimensions $\le 6$ (and perhaps in general) that admit
symplectic structures.
\end{abstract}

\maketitle

\section{Introduction}

Symplectic structures, which are closed, nondegenerate 2-forms $\omega$
on an even-dimensional manifold $X$, can be thought of as skew-symmetric 
analogs of constant curvature Riemannian metrics. 
The nondegeneracy condition (that each nonzero tangent vector pairs 
nontrivially with some vector) is the same in each case, and closure 
($d\omega=0$) corresponds to a constant curvature condition in that it is 
a differential equation guaranteeing that all such structures of a given 
dimension are locally identical. 
The study of constant curvature manifolds reduces, through covering space 
theory, to that of discrete isometry groups of Euclidean, hyperbolic and 
spherical space, so it is natural to ask whether symplectic structures 
also have some sort of topological characterization. 
Gromov \cite{Gr}, \cite[Theorem 7.34]{McS}, showed that an open 
manifold $X$ admits a symplectic structure for each choice of 
almost-complex structure (up to homotopy) and class 
$[\omega]\in H_{\dR}^2(X)$, reducing the existence question for symplectic 
structures on open manifolds to that of almost-complex structures 
(homotopy theory for the tangent bundle). 
However, the case of closed manifolds is much more difficult. 
For example, there exist homeomorphic pairs of smooth 4-manifolds with
isomorphic tangent bundles, such 
that one admits symplectic structures and the other does not \cite{Ta}. 
(See also \cite{GS}, \cite{K}.)
In the present paper, we propose a solution to this problem, by 
introducing a topological structure called a hyperpencil, which we show 
determines a symplectic structure (up to isotopy). 
In dimensions $\le6$, and conjecturally in general, a closed manifold 
admits a symplectic structure if and only if it admits a hyperpencil, 
and a dense subset of all symplectic structures (up to isotopy and scale) 
can be described as a quotient of the set of deformation classes of 
hyperpencils. 

The prototype for hyperpencils comes from algebraic geometry 
(Example~\ref{example:2.6}(a)). 
If $X\subset \CP^N$ is a smooth $n$-dimensional algebraic variety, we 
obtain a {\em linear $k$-system\/} $f:X-B\to \CP^k$ on $X$ by 
intersecting $X$ with a transverse linear subspace $A$ of codimension 
$k+1$, setting $B= X\cap A$, and defining $f$ to be the restriction 
of a projectivized linear surjection $\CP^N-A\to \CP^k$. 
Thus, any algebraic variety inherits a canonical deformation class of 
linear $k$-systems from its embedding in $\CP^N$. 
One can use the resulting local structure on $X$ to formulate a definition 
of linear systems in the category of smooth manifolds. 
(\cite{G3} studies the general case of this.) 
When $k$ equals 1 and $n$, respectively, generic prototypes yield Lefschetz
pencils and (singular) branched coverings, both of which have been extensively
studied by topologists.
Lefschetz used Lefschetz pencils to study the topology of algebraic 
varieties (e.g., \cite{L}), and in recent decades these structures have 
also arisen in 4-manifold theory (e.g., \cite{GS}). 
We wish to use linear systems to construct symplectic structures. 
To obtain the strongest theorem, we wish to use the weakest possible 
hypotheses. 
This suggests using the smallest possible value for $k$, since a linear 
$k$-system generates linear $\ell$-systems for $\ell\le k$
(assuming the associated almost-complex structures behave reasonably 
as in the algebraic case) by composition with a generic projection 
$\CP^k - A' \to \CP^\ell$. 
However, the fibers of a linear $k$-system have (real) dimension $2(n-k)$; 
when $n-k>1$ it is already difficult to know when the fibers admit 
symplectic  structures. 
Thus, the optimal case seems to be when $k=n-1$, when generic fibers are 
oriented surfaces, so each has a unique symplectic form (i.e., area form 
in this dimension) up to isotopy and a constant scale factor. 
Hyperpencils (Definition~\ref{def:2.4}) are a type of linear 
$(n-1)$-system derived from the algebraic prototype. 
We have aimed for the weakest possible hypotheses guaranteeing the existence 
of symplectic structures, allowing the ugliest possible local behavior. 
It seems likely that additional constraints should be added for other 
purposes; for example, it may be possible to deform any hyperpencil 
into a much nicer ``generic'' form. 

Our Main Theorem~\ref{thm:2.11} can be paraphrased as follows: 

\begin{thm}\label{thm:1.1} 
Let $X$ be a smooth, closed, oriented manifold. 
\begin{itemize}
\item[(a)] For any hyperpencil on $X$, the space of suitably compatible 
almost-complex structures $J$ is nonempty and contractible.
\item[(b)] Every such $J$ is tamed by a symplectic form on $X$ 
realizing a certain cohomology class determined by $f$, and the 
space of such forms is contractible (for $J$ fixed or suitably varying).
\item[(c)] There is a canonical 
map $\Omega :\P(X) \to \cS(X)$, where $\P(X)$ is the set of deformation 
classes of hyperpencils on $X$, and $\cS(X)$ is the set of all 
isotopy classes of symplectic forms on $X$.
\end{itemize}
\end{thm}

We first consider (c). 
Deformations of hyperpencils are defined in Definition~\ref{def:2.7}. 
Symplectic forms $\omega_0$ and $\omega_1$ are {\em isotopic\/} if 
there is a diffeomorphism $\varphi :X\to X$ isotopic to $\id_X$ with 
$\varphi^* \omega_1=\omega_0$. 
By Moser's Theorem \cite{M}, this is equivalent to the existence 
of a deformation (smooth family of symplectic forms $\omega_s$, $0\le
s\le1$) for which $[\omega_s] \in H_{\dR}^2 (X)$ is constant. 
Theorem~\ref{thm:2.11}(c)  characterizes  the 
symplectic forms associated by $\Omega$ to a given hyperpencil,  
using the intermediate structure $J$ 
(see Definition~\ref{def:2.9} and Lemma~\ref{lem:2.10}). 
For example, it is easy to check that the standard K\"ahler form on an 
algebraic variety is associated in this manner to the deformation 
class of hyperpencils determined by its embedding in $\CP^N$. 
A more expository discussion of Theorem~\ref{thm:2.11} (in a slightly 
earlier form) appears in \cite{G2}. 
The original form of the theorem, that a 4-manifold with a Lefschetz pencil
admits a symplectic structure,was first proved by the author in 1990, but 
remained unpublished (due to the 
emergence of a more  direct way of constructing unusual symplectic 
4-manifolds \cite{G1}) until its expository appearance as
\cite[Theorem 10.2.18 and Corollary 10.2.23]{GS}. 

The main motivation for Theorem~\ref{thm:1.1} is the use of (c) in 
characterizing symplectic manifolds. 
The symplectic forms produced by the theorem are integral (i.e., with
cohomology class in the 
image of $H^2 (X;\zed)\to H_{\dR}^2(X)$). 
Donaldson \cite{D} has proven that any integral symplectic manifold
(up to scale)
admits an associated Lefschetz pencil, and Auroux \cite{A1} has obtained
a similar result for linear 2-systems. 
These results imply the $n=2,3$ cases, respectively, of the following 
conjecture (which is trivially true for $n\le1$):

\begin{conj}\label{conj:1.2} 
Let $\omega$ be any integral symplectic form on a closed $2n$-manifold $X$. 
Then for any sufficiently large integer $m$, the isotopy class of 
$m\omega$ lies in the image of $\Omega$.
\end{conj} 

\noindent
The conjecture is still open for $n\ge4$, motivating our attempt at the 
weakest possible definition of hyperpencils. 
However, Auroux has made some technical progress on the problem \cite{A2}. 
One would ultimately expect hyperpencils arising from Donaldson-Auroux 
theory to have much nicer local properties than arbitrary hyperpencils, 
for example explicit holomorphic local models at the critical points. 
The conjecture leads to characterization of symplectic manifolds as follows.
Up to scale, every rational cohomology class is integral, 
and the subspace $\cS_{\que}(X) \subset \cS(X)$ of symplectic forms with
$[\omega]$ rational is dense (since nondegeneracy is an 
open condition). 
Thus, we may define $\tilde \Omega :\P(X) \times \que_+ \to \cS(X)$ so that 
$\tilde\Omega (\varphi, q)$ is obtained by rescaling $\Omega(\varphi)$ to make its 
cohomology class $q$ times a primitive integral class, and conclude: 

\begin{prop}\label{prop:1.3} 
If Conjecture~\ref{conj:1.2} holds for $X$, then the image of the 
canonical map $\tilde\Omega:\P(X) \times \que_+ \to \cS(X)$ is the dense 
subset $\cS_{\que} (X) \subset\cS(X)$.
\end{prop}

\begin{cor}\label{cor:1.4} 
In dimensions where Conjecture~\ref{conj:1.2} holds (e.g., dimensions 
$\le6$), a closed manifold admits a symplectic structure if and only if 
it admits a hyperpencil. 
A closed 4-manifold admits a symplectic structure if and only if it admits 
a Lefschetz pencil with $B\ne\emptyset$. 
\end{cor}

\noindent 
The last statement of the corollary follows from 
Theorem \ref{thm:1.1}(c) and \cite{D}, if we restrict to Lefschetz
pencils such that each irreducible component of each singular 
fiber intersects $B$ (since these are hyperpencils); the stated
version is proved directly in 
\cite[Theorem 10.2.28]{GS}.
(There, the condition $B\neq \emptyset$ is contained in the definition of
Lefschetz pencils.) 
In dimensions where Conjecture~\ref{conj:1.2} holds, we have now 
topologically characterized manifolds admitting symplectic structures. 
{From} there, to topologically determine the dense subset $\cS_{\que}(X)
\subset \cS(X)$, it suffices to understand the following: 

\begin{conj}\label{conj:1.5} 
The fibers of $\tilde\Omega$ (or equivalently, of the map $\overline{\Omega}:
\P(X) \to \cS_\que(X)/\que_+$ determined by $\Omega$) are specified by a 
topologically defined equivalence relation on $\P(X)$. 
\end{conj}

\noindent
This may be easier to prove for a stronger definition of hyperpencils. 
The main evidence for Conjecture~\ref{conj:1.5} is that the theorems 
of Donaldson and Auroux come with uniqueness statements up to a notion 
of stabilization, which multiplies the cohomology classes by large integers.
While this stabilization comes from analytical considerations on special 
families of linear systems, one might hope to topologically define 
stabilization maps in general, $\sigma_k :\P(X)\to \P(X)$, $k\in\zed_+$, 
with $\sigma_1 =\id_{\P(X)}$, $\sigma_k \circ\sigma_\ell = \sigma_{k\ell}$ 
and $\Omega \circ\sigma_k = k\Omega$, and realize the equivalence 
relation in Conjecture~\ref{conj:1.5} by the definition $\varphi\sim \psi$ if 
and only if $\sigma_k(\varphi) = \sigma_\ell(\psi)$ for some $k,\ell\in\zed_+$. 
However, these stabilizations already seem complicated in dimension~4. 
(For $\sigma_2$, see \cite{AK}.) 

Our main tool for constructing  symplectic structures is a method 
originally used by Thurston \cite{T} in the context of surface bundles, to use 
a symplectic structure on the base to construct one on the total space. 
This method has been generalized to bundles with higher dimensional 
fibers (e.g., \cite[Theorem 6.3]{McS}) and to bundles with complex 
quadratic singularities \cite[Theorem 10.2.18]{GS}, but we show 
(Theorem~\ref{thm:3.1}) that the method works for maps that may be very 
different from bundle projections. 
For example, it suffices to have a map that is $J$-holomorphic for 
suitable almost-complex structures, with suitable data in a neighborhood 
of each point preimage. 
(An \emph{almost-complex structure} on a manifold $X$ is a complex structure 
on its tangent bundle $TX$, or equivalently a bundle map $J:TX\to TX$ 
covering $\id_X$ with $J\circ J= -\id_{TX}$, which we should interpret 
as multiplication by $i$.
A map is then $J$-\emph{holomorphic} if its derivative is complex linear.
Complex structures on other vector bundles may be interpreted similarly.) 
We apply this method to a hyperpencil $f:X-B\to \CP^{n-1}$, starting 
from the standard symplectic form $\ostd$ on $\CP^{n-1}$. 
To do this, we need a suitable almost-complex structure on $X$, and to 
prove uniqueness up to isotopy we must be able to find a 1-parameter 
family $J_s$ connecting any two such almost-complex structures. 
Thus, we need various lemmas for splicing together locally defined 
almost-complex structures. 
These are compiled into Lemma~\ref{lem:3.2}, whose proof comprises 
most of Section~4. 
To emphasize that the choice of almost-complex structure does not 
crucially affect the resulting symplectic forms, we define hyperpencils 
using only locally defined almost-complex structures and prove that 
the relevant space of global almost-complex structures is 
nonempty and contractible
(Theorem~\ref{thm:1.1}(a)). 
For convenience, and to emphasize the topological nature of the 
hypotheses, we always work with $C^0$ almost-complex structures. 
Thus, our spaces of almost-complex structures will always be given the 
$C^0$-topology (or for noncompact $X$, its natural generalization, the 
compact-open topology). 
In contrast, we have much more flexibility in topologizing spaces of 
symplectic forms. 
For example, contractibility in Theorem~\ref{thm:1.1}(b) holds for all 
$C^k$-spaces of forms and Sobolev spaces in between. 
(See Theorem~\ref{thm:2.11}(b).)
Our method for constructing symplectic structures has other applications 
besides Theorem~\ref{thm:1.1}. 
We study high-dimensional Lefschetz pencils and other linear systems in 
\cite{G3}, and locally holomorphic maps with 2-dimensional fibers in \cite{G4}.

Throughout the paper, orientations are crucial. 
If $V$ is a $2n$-dimensional real vector space, any nondegenerate, 
skew-symmetric, bilinear form $\omega$ on $V$ induces an orientation, 
since its top exterior power is a volume form. 
A (linear) complex structure $J$ on $V$ induces an orientation obtained, as 
usual, from any complex isomorphism $(V,J) \cong \complex^n = \complex 
\times \cdots\times \complex$ in the product orientation (where $(1,i)$ 
is a positively  oriented real basis for $\complex$). 
If $V$ is given to be an oriented vector space, we only consider 
(unless otherwise specified) complex structures and nondegenerate 2-forms 
inducing the given orientation on $V$. 
For example, almost-complex structures and symplectic forms on oriented 
manifolds implicitly induce the given orientation. 
We let $\ostd$ denote the standard symplectic form on $\CP^k$, normalized 
so that $\int_{\CP^1} \ostd =1$, so $[\ostd] \in H_{\dR}^2 (\CP^k)$ 
is the hyperplane class $h$, Poincar\'e dual to $[\CP^{k-1}] \in H_{2k-2}
(\CP^k)$. 

The author wishes to acknowledge Denis Auroux for helpful discussions. 

\section{The main theorem}
To define hyperpencils, we need some preliminary definitions.
We begin by generalizing some standard terminology for relating
symplectic and complex structures.

\begin{de} \label{def:2.1}
Let $T:V \to W$ be a linear transformation between finite-dimensional
real vector spaces, and let $\omega$ be a skew-symmetric bilinear form on $W$.  
A linear  complex structure $J:V \to V$ will be called 
$(\omega,T)$-{\em tame\/} if $T^*\omega(v, Jv) > 0$ for all $v \in V
- \ker T$.  If, in addition, $T^*\omega$ is $J$-invariant (i.e., 
$T^*\omega(Jv, Jw) = T^*\omega(v,w)$ for all $v,w \in V$), we
will call $J$ $(\omega, T)$-{\em compatible\/}.  For a $C^1$-map $f:X
\to Y$ between manifolds, with a 2-form $\omega$ on $Y$, an
almost-complex structure $J$ on $X$ will be called
$(\omega,f)$-{\em tame\/} (resp. $(\omega, f)$-{\em compatible\/}) if it
is $(\omega,df_x)$-{\em tame\/} (resp. $(\omega,df_x)$-{\em compatible\/}) 
for each $x \in X$.  If $T = \id_V$ or $f = \id_X$, we will shorten
the terminology to $\omega$-{\em tame\/} and $\omega$-{\em compatible\/}. 
\end{de}

\noindent 
The last sentence of the definition is standard terminology.  
If $\omega$ tames some $J$  (i.e., $T=\id_V$ and 
$J$ is $\omega$-tame), then $\omega$ is
obviously nondegenerate, so a closed, taming 2-form $\omega$ is automatically
symplectic.  
An $\omega$-tame $J$ induces the same orientation as
$\omega$. (For example, homotope $\omega$ through taming forms to a
compatible one.)  
A key advantage of $\omega$-taming over $\omega$-compatibility is that 
the former is an open condition on manifolds.  
In fact, the $\omega$-taming condition is satisfied
provided that it holds on the unit sphere bundle in $TX$, so it is
preserved under $\varepsilon$-small perturbations of $J$ and $\omega$
when $X$ is compact.
If $J$ is $(\omega, T)$-tame then $\ker T$ is a
$J$-complex subspace of $V$, so $\Im T \subset W$ inherits an
$\omega$-tame complex structure $T_*J$ making $T$ complex linear
($T \circ J = T_*J \circ T$).  
This will be $\omega$-compatible if
and only if $J$ is $(\omega, T)$-compatible.  
For an
$(\omega,f)$-tame almost-complex structure, preimages of regular
values of $f$ will be $J$-holomorphic submanifolds (i.e., $J$
preserves their tangent spaces), and the complex
structures induced on the fibers of $\Im df \subset f^*TY\to X$ 
will be denoted $f_*J$.  
Both the taming and compatibility
conditions are preserved under taking convex combinations $\sum t_i
\omega_i$  (all $t_i \geq 0, \sum t_i = 1$) for fixed $f,J$. 
An almost-complex symplectic manifold $(Y,J,\omega)$ is called 
{\em almost-K\"ahler\/} if $J$ is $\omega$-compatible, and 
{\em K\"ahler\/} if,  in addition, $(Y,J)$ is a complex manifold. 
In either case, if $f:X\to Y$ is $J$-holomorphic for some 
almost-complex structure on $X$, this structure is $(\omega,f)$-compatible. 

To prove uniqueness of symplectic forms induced by hyperpencils, we
need a technical condition for critical points.  Suppose $E,F \to X$
are real (finite dimensional) vector bundles over a metrizable 
topological space,
and $T:E \to F$ is a (continuous) section of the bundle $\Hom(E,F)$. 
In our main application, these will be induced by a $C^1$-map $f:X \to Y$
between manifolds, with $T= df:TX \to f^*TY$. Motivated by
this example, we call a point $x \in X$ {\em regular\/} if $T_x:E_x \to
F_x$ is onto and {\em critical\/} otherwise.  
Let $P \subset E$ be the
closure $\cl( \bigcup \ker T_x)$, where $x$ varies over all the regular
points of $T$ in $X$, and let $P_x = P \cap E_x$.  Thus, $P_x = \ker
T_x$ if $x$ is regular, and otherwise $P_x \subset \ker T_x$ consists
of limits of sequences of vectors annihilated by $T$ at regular
points. 

\begin{de} \label{def:2.2}
A point $x \in X$ is {\em wrapped\/} if $\Span P_x$ has (real)
codimension at most 2 in $\ker T_x$.
\end{de}

\begin{prop} \label{prop:2.3}
Suppose that in a neighborhood of a critical point $x \in X$, $T$ is
given by $df$, for some holomorphic map $f:U \to
\complex^{n-1}$ with $U$ open in $\complex^n$.  If each fiber
$f^{-1}(y)$ intersects the critical set $K$ of $f$ in at most a finite
set, then $x$ is wrapped.  
In fact, $P_x = \ker T_x$.  
\end{prop}

\noindent
This proposition will show that our hypothesis of wrapped critical
points is broad enough to be useful.  Note, however, that the
proposition becomes false without the finiteness hypothesis, e.g.,
$n=3$, $f(x,y,z)=(x^2, y^2)$ at $(0,0,0)$.  
(For $n=2$, $P_x$ equals $\ker T_x$ unless $f$ is constant or $x$ is a smooth
point of $f$ with multiplicity $>1$; 
cf. \cite[proof of Proposition 1.3]{G4}.)
Similarly, $P_x$ may not
equal $\ker T_x$ if we pass from the holomorphic setting to $C^\infty$. 
(The $C^\infty$-map $f:\real^2 \to \real$ given by
$f(x,y) = y^3 + e^{-1/x^2}y$ has a unique critical point at $(0,0)$,
so $df_{(0,0)} = 0$, but $P_{(0,0)}$ is the $x$-axis.)

\begin{proof}
For $\ell \geq 2$, let $K_\ell \subset K$ be the set of $z \in U$ for
which $\ker df_z$ has complex dimension $\geq \ell$.  Thus, 
$K = K_2 \supset K_3 \supset \cdots \supset K_{n+1} = \emptyset$.  We
begin by showing that each $K_\ell$ is an analytic variety of complex
codimension $\geq \ell$ in $U$.  Analyticity follows immediately from
the description of $K_\ell$ as the set of $z \in U$ for which every
$(n-\ell+1) \times (n-\ell+1)$ submatrix of $df_z$ has
determinant zero.  For $z$ in the top stratum $W$ of $K_\ell$, let
$Q_z = T_z W \cap \ker df_z$. If $\codim_\complex K_\ell <
\ell$, then $Q_z$ has nonzero dimension for all $z \in W$.  
Choose some
$z_0 \in W$ minimizing this dimension.  
To see that $Q$ is a smooth
distribution on $W$ near $z_0$, choose a projection $\pi$ of
$\complex^{n-1}$ whose restriction to $df_{z_0}(T_{z_0}W)$ is
an isomorphism, and note that $\pi \circ f|W$ is a submersion at $z_0$
with $\ker d(\pi \circ f|W)_z$ containing $\ker d(f|W)_z = Q_z$; 
these latter spaces are then equal near $z_0$ by minimality of $\dim Q_{z_0}$.  
Now choose a smooth, nonzero vector field in $Q$ near $z_0$.  
By integrating, we obtain a curve in $W \subset K$ whose image
under $f$ is a point $y$, contradicting finiteness of $f^{-1}(y) \cap K$. 

Now observe that the subset $V=\bigcup_{z \in U} \ker df_z \subset
TU = U \times \complex^n$ is an analytic variety with complex
dimension $\geq n + 1$ everywhere, since it is cut out by the system
of $n-1$ equations $df_z(v)=0$ in $(z,v)$.  
For each $\ell \geq
2$, the subset $V_\ell = V \cap ((K_\ell -K_{\ell+1}) \times
\complex^n)$ is a complex $\ell$-plane bundle over $K_\ell -
K_{\ell+1}$, and the latter has codimension $\geq \ell$ in $U$, so
$\dim_\complex V_\ell \leq n < \dim_\complex V$.  
Thus, $V - K \times \complex^n = V - \bigcup_{\ell=2}^n V_\ell$ 
has closure $V$ in $U\times\complex^n$.
The proposition follows immediately.
\end{proof}

\begin{de} \label{def:2.4}
A {\em hyperpencil\/} on a smooth, closed, oriented, $2n$-manifold $X$
is a (necessarily finite) subset $B \subset X$ called the {\em base
locus\/} and a smooth map $f:X-B \to \complex  \PP^{n-1}$ such that 
\begin{itemize}
\item[(1)] each $b \in B$ is mapped to $0 \in \complex^n$ by an
  orientation-preserving local coordinate chart in which $f$ is given
  by projectivization $\complex^n-\{0\} \to \complex \PP^{n-1}$,
\item[(2)] each critical point of $f$ is wrapped and has a
  neighborhood with a continuous $(\omega_{\std}, f)$-compatible
  almost-complex structure, and
\item[(3)] each fiber $F_y = \cl f^{-1}(y) \subset X$ contains only
  finitely many critical points of $f$, and each component of each
  $F_y - \{\text{critical points}\}$ intersects $B$.
\end{itemize}
\end{de}

\begin{Remarks} \label{rem:2.5}
(a) The results in this paper are all trivially true for $n\le1$ (after 
Moser \cite{M}), so we will assume $n\ge 2$ whenever convenient.
For potential applications such as sub-hyperpencils, where it may be 
convenient to allow $n\le1$, we specify the required conventions: 
For $n=0$, $B$ equals $X$, and a symplectic form on a 0-manifold is its 
unique positive orientation. 
For $n=1$, we require $B\subset X$ to be finite, and its Poincar\'e 
dual in $H_{\dR}^2(X)$ is the class $c_f$ used in Theorem~2.11. 
\smallskip

\noindent (b) By Condition (1), the fibers $F_y$ of a hyperpencil 
are complex lines near each $b \in B$.  
Thus, $F_y = f^{-1}(y) \cup B$ and $F_y \cap F_z = B$ for $y
\neq z$, so $B$ is closed and discrete, hence finite.  
Each $F_y$ is an
oriented surface (in the preimage orientation) except possibly at  
finitely many singularities, where $F_y$ intersects the critical set of 
$f$ in $X-B$.  
\smallskip

\noindent (c) 
For $n \le 2$ the hypothesis of wrapped critical points is trivially true.  
(The proof of Main Theorem~2.11 shows that the
regular points of $f$ are dense in $X-B$, so each $\Span_\real P_x \subset
T_x(X-B)$ is a nontrivial complex subspace.)  
For $n=3$, this
hypothesis can be eliminated if we assume that at each point in the
closure of the set of 
unwrapped critical points, the given local
almost-complex structure makes $f$ 
$J$-holomorphic for some continuous, $\omega_{\std}$-tame local
complex structure on the bundle $f^*T\complex\PP^{n-1}$. 
(In fact, an even weaker hypothesis guarantees the Main Theorem when $n=3$, 
namely $(\omega_{\std},df)$-extendability as used in 
Addendum~\ref{addendum} 
with $C=D=\emptyset$ and $E= TX$.)
For arbitrary $n$, the hypothesis of wrapped critical points 
can be dropped in the presence of a global $\omega_{\std}$-compatible 
complex structure (standard near $B$) on  $f^*T\complex\PP^{n-1}$ 
making $f$  $J$-holomorphic for each local almost-complex
structure on $X-B$, but the resulting isotopy class of symplectic
forms could then conceivably depend on the choice of this
structure on $f^* T\complex \PP^{n-1}$. 
\smallskip

\noindent 
(d) Throughout the article, we use continuous, rather than smooth,
almost-complex structures.  This is both for convenience (avoiding
awkward and unnecessary proofs of smoothness) and to emphasize that the
purpose of the local almost-complex structures is topological rather
than analytical, controlling monodromy around the critical values.
For a Lefschetz pencil on a 4-manifold, for example, the monodromy
consists of right-handed Dehn twists.  Allowing the opposite
handedness violates the hypothesis of (compatibly oriented) local
almost-complex structures, and results in manifolds having no
symplectic structure. 
\smallskip

\noindent
(e) It is an open question whether $(\omega_{\std},f)$-compatibility 
can be replaced by $(\omega_{\std},f)$-taming.  This could be
done throughout the paper 
if Question \ref{ques:4.3} had an affirmative answer, and can
also be done in the situation at the end of (c) above (arbitrary $n$).
\end{Remarks}

\begin{examples} \label{example:2.6}
(a) Any smooth algebraic variety $X \subset \complex \PP^N$ admits a
hyperpencil.  Simply pick a linear subspace $A \approx \complex
\PP^{N-n}$ $(n=\dim_\complex X)$ transverse to $X$ in $\complex
\PP^N$, let $B=X \cap A$ and let $f$ be the restriction of the
holomorphic projection $\complex \PP^N - A \to \complex
\PP^{n-1}$. Condition (1) follows from transversality. For Condition
(3), note that each fiber $F_y$ is an algebraic curve in some linear
subspace of $\complex \PP^N$ containing $A$ as a codimension-1
subspace.  
Thus, each irreducible component of $F_y$ intersects $A$,
and so has a 0-dimensional intersection with the critical set 
(since Condition (1) locates some regular points).  
Now the critical points of $f$ are wrapped by
Proposition~\ref{prop:2.3}, and the obvious holomorphic structure on
$X$ completes Condition (2). 

\noindent
(b) On 4-manifolds, {\em Lefschetz pencils\/} are inspired by the above
Example (a) with $n=2$.  
Their definition is obtained from Definition~\ref{def:2.4} 
$(n=2)$ by replacing Conditions (2) and (3) with the
condition that $f$ should  locally be a complex Morse function, i.e., 
modeled by $f(z_1, z_2) = z_1^2 + z_2^2$ at each critical point.  (This
holds in (a) for a generic choice of $A$.) See \cite{GS} for a recent
survey on Lefschetz pencils. 
Note that Condition (2) is automatically
satisfied by Lefschetz pencils, as is the finiteness part of Condition~(3). 
Lefschetz pencils need not satisfy the other part of Condition (3), however. 
This condition on $F_y \cap B$, 
which guarantees that a (4-dimensional) Lefschetz pencil is a hyperpencil, 
is actually not necessary for constructing symplectic
structures on Lefschetz pencils. 
(See \cite[Theorem~10.2.18 and 
Corollary~10.2.23]{GS} for details.) It suffices to know either
that $B \neq \emptyset$ or that the fibers are nontrivial in $H_2(X;
\real)$.  
(In fact, the only counterexamples are torus bundles, namely
$L(p, 1) \times S^1 \to S^2$, and their blowups; see
\cite[Remark~10.2.22(a)]{GS}.)  
However, without Condition (3), one loses control of the
cohomology class of the resulting form, and with it one loses
uniqueness of the isotopy class.
\end{examples}

To define an appropriate equivalence relation among hyperpencils on
$X$, we again work by analogy with algebraic geometry.  
We begin by
organizing hyperpencils into families over a parameter space $S$.
Roughly, these families are given by bundles over $S$ whose fibers
have continuously varying hyperpencil structures.

\begin{de} \label{def:2.7}
A {\em family of hyperpencils\/} parametrized by a topological space
$S$ consists of a pair of fiber bundles $\pi_X:X \to S$, $\pi_Y:Y \to
S$, a subset $B \subset X$ and a continuous, fiberwise
smooth map $f:X-B \to Y$ covering $\id_S$, subject to the following:  
The fibers $X_s$ of $X$, $s \in S$,  are all diffeomorphic to a
fixed, closed, oriented $2n$-manifold (so the structure group consists
of orientation-preserving diffeomorphisms of the fiber in the
$C^\infty$-topology), and the fibers
of $Y$ are diffeomorphic to $\complex \PP^{n-1}$ with structure group
$\PP \rU(n)$ acting in the usual way.  The map $\pi_X|B:B \to S$ is a
(necessarily finite) covering map. 
In addition:
\begin{itemize}
\item[(1)] $B$ has a neighborhood $V \subset X$ on which $\pi_X|V:V
  \to S$ lifts to $\tilde{\pi}:V \to B$, and $\tilde{\pi}$ is given the
  structure of a $\rU(n)$-vector bundle (with zero section $B$ and fibers
  oriented compatibly with those of $\pi_X$).  The map $f|V-B$ is
  projectivization on each fiber of $\tilde{\pi}$.
\item[(2)] The map $df:T^v(X-B) \to f^*T^v Y$ (where $T^v$
  denotes the bundle of tangent spaces to the fibers of $\pi_X$ and
  $\pi_Y$) is continuous, and  
each critical point of $df$ has a neighborhood in $X-B$ over which
  $T^vX$ has an $(\omega_{\std}, df)$-compatible complex
  vector bundle structure.
\item[(3)] For each fiber $X_s$ of $\pi_X$, each critical point of
  $f|X_s - B$ is wrapped (in $X_s$), and Condition (3) of 
Definition~\ref{def:2.4} is satisfied.
\end{itemize}
A {\em deformation\/}  of hyperpencils is a family parametrized by
$I = [0,1]$ with $X=I \times X_0$. 
\end{de}

It is easily verified that a family of hyperpencils parametrized by a
1-point space is the same as a hyperpencil (together with a fixed choice 
of the charts in Definition~\ref{def:2.4}(1) up to $\rU(n)$ action).  
If $\varphi:S \to S'$ is
continuous, then a family of hyperpencils parametrized by $S'$ pulls
back to one parametrized by $S$.  
For example, any parametrized
family of hyperpencils restricts to a hyperpencil on each $X_s$, 
or to a family parametrized by any subspace of $S$. 
Now, we easily obtain an
equivalence relation by calling two hyperpencils on a fixed manifold
{\em deformation equivalent\/} if they are realized as $X_0$ and $X_1$
for some deformation.  
(The only technicality is that in transitivity, the middle hyperpencil 
may inherit two different sets of the charts in Definition~\ref{def:2.4}(1).
However, these charts can easily be changed, by triviality of bundles over $I$.
Closer inspection also shows that we can find a continuous family 
interpolating between any two such charts; cf.\ proof of Lemma~2.10.) 
For a parametrized family with $S$ 
path connected, any $X_s$ and $X_t$ can be identified so 
that their hyperpencils are deformation equivalent. 
We can
also assume the deformation is constant near each endpoint of $I$.
If, in addition, $X=S \times X_s$ is the trivial bundle, we obtain a
family of deformation equivalent hyperpencils on the fixed manifold $X_s$. 

\begin{examples} \label{example:2.8}
(a) To construct families of hyperpencils as in Example \ref{example:2.6}(a), 
let $G$ denote the complex Grassmann manifold of codimension-$n$
linear subspaces of $\complex \PP^N$, and let $\gamma \subset G \times
\complex \PP^N$ be the tautological bundle whose fiber over $A \in G$
is $A \subset \complex \PP^N$.  If $Y \to G$ is the bundle whose
fiber over $A$ is the $(n-1)$-plane  in $\complex \PP^N$ with
maximal distance from $A$, we obtain a canonical holomorphic map $f:G
\times \complex \PP^N - \gamma \to Y$ covering $\id_G$, induced by
linear projections on $\complex^{N+1}$. Given an $n$-dimensional
smooth algebraic variety $X_0$ in $\complex \PP^N$, let $S \subset G$
be the open set of subspaces transverse to $X$.  Restricting $f$ to a
map $S \times X_0 - \gamma \to Y|S$, we obtain a parametrized family of
hyperpencils, consisting of all hyperpencils on $X_0$ obtained by
Example \ref{example:2.6}(a).  (The given holomorphic structure on $S
\times X_0$ satisfies Condition (2) above.)  For example, when $n=2$,
generic members of the family will be Lefschetz pencils, but there
will typically also be parameter values for which quadratic critical
points coalesce into those of higher degree.  
In general, the space
$S$ is path connected (since it is obtained from $G$ by removing a
subvariety of positive complex codimension), so we conclude that all
hyperpencils on $X_0$ obtained by Example~\ref{example:2.6}(a) are
deformation equivalent (for a fixed embedding $X_0 \subset \complex
\PP^N$). 
\smallskip

\noindent (b) For a family in which $X$ is a nontrivial bundle, note
that the space of all hypersurfaces of a fixed degree in $\complex
\PP^N$ is parametrized by some $\complex \PP^M$.  
Let $S \subset
G \times \complex \PP^M$ denote the path-connected subset of pairs
$(A, t)$ such that the variety $X_t$ is nonsingular and transverse to $A$.  
The construction of (a) above generalizes immediately to produce
a parametrized family consisting of all hyperpencils as above 
on all nonsingular hypersurfaces of a fixed degree.  
It follows that any two nonsingular
hypersurfaces of the same degree in $\complex \PP^N$ are diffeomorphic
in such a way that the canonical families of hyperpencils are
deformation equivalent.  We also see that the canonical deformation
class of hyperpencils on a fixed hypersurface is invariant under
self-diffeomorphisms induced by monodromy of the bundle of all
nonsingular hypersurfaces.
\end{examples} 

We relate hyperpencils $f:X- B\to \complex \PP^{n-1}$ 
to symplectic structures $\omega$ on $X$ via the existence of local 
almost-complex structures that are simultaneously $\omega$-tame and 
$(\omega_{\std},f)$-compatible: 

\begin{de} \label{def:2.9}
Let $\omega$ be a continuous 2-form on an oriented manifold $X$.  
A hyperpencil
$f:X - B \to \complex \PP^{n-1}$ is $\omega$-{\em tame\/} if $X$
is covered by open sets $W_\alpha$ with continuous, $\omega$-tame
almost-complex structures $J_\alpha$ such that $J_\alpha|W_\alpha - B$ 
is $(\omega_{\std}, f)$-compatible, and the structures $f_*
J_\alpha$ on $\Im df \subset f^*T \complex \PP^{n-1}|
W_\alpha - B$ all agree where their domains overlap. Similarly, a
family $f:X-B \to Y$ of hyperpencils parametrized by $S$ is
$\omega$-{\em tame\/} for a continuous family $\omega$ of 2-forms on
the fibers of $\pi_X$, if $X$ is covered by $(W_\alpha, J_\alpha)$ as
above.  (Here $W_\alpha$ is open in $X$, $J_\alpha$ is a complex
vector bundle structure on $T^v X|W_\alpha$, and the structures
$f_*J_\alpha$ fit together on $\Im df \subset f^*T^v Y$.)
\end{de}

\begin{lem} \label{lem:2.10}
A hyperpencil $f$ is $\omega$-tame if and only if there is a global,
continuous, $\omega$-tame almost-complex structure $J$ on $X$ with
$J|X - B$ $(\omega_{\std}, f)$-compatible.
The structure $J$ can be chosen to agree
near $B$ with the local complex structures induced by any preassigned charts 
as in Definition~\ref{def:2.4}(1), 
and so that $f_*J$ agrees with the structures $f_*J_\alpha$ given 
by Definition~\ref{def:2.9} if  these are standard near $B$.  
The corresponding statements also hold for a 
family of hyperpencils parametrized by a metrizable space,
where $J$ is a complex vector bundle structure on $T^vX$ agreeing on a
neighborhood of $B$ with the fiberwise complex structure determined by
the vector bundle in Definition~\ref{def:2.7}(1).
\end{lem}

\noindent
In other words, the local structures $J_\alpha$ of Definition~\ref{def:2.9} 
can be patched together to form a global structure. 
The ``if'' direction of this lemma is obvious; the other is proved in
Section 4.

Any hyperpencil $f:X - B \to \complex \PP^{n-1}$ determines a 
class $c_f = f^*h \in H^2(X; \integer) \cong H^2(X - B; \integer)$,
where $h$ is the hyperplane class, Poincar\'{e} dual to $[\complex
\PP^{n-2}] \in H_{2n-4}(\complex \PP^{n-1}; \integer)$.  
This class in $X$ is invariant under deformations of $f$.  
When $X \subset
\complex \PP^N$ is a smooth algebraic variety with its canonical
deformation class of hyperpencils, $c_f$ is the restriction of the
hyperplane class of $\complex \PP^N$, and in general for $n=2$, $c_f$
is Poincar\'{e} dual to any fiber $F_y$.  
We also use $c_f$ to
denote the corresponding class in $H^2_{\dR}(X)$; 
this is the cohomology class of the symplectic forms associated to $f$. 

We are now ready to state the Main Theorem.		

\begin{thm} \label{thm:2.11}            
Let $X$ be a smooth, closed, oriented $2n$-manifold.
\begin{itemize}
\item[a)] For any hyperpencil $f:X - B \to \complex \PP^{n-1}$ there
  is a continuous almost-complex structure $J$ on $X$ that is
  $(\omega_{\std}, f)$-compatible on $X - B$, and agrees near
  $B$ with the complex structure given there by Definition~\ref{def:2.4}(1) 
(for any fixed choice of charts). 
The space of all such structures $J$ in the
  $C^0$-topology (for fixed $f$ and charts near $B$) 
is contractible, as is the space of
  all $J$ on $X$ that are $(\omega_{\std}, f)$-compatible on
  $X - B$ (without the constraint near $B$).
\item[b)] For any $J$ as in the first sentence of (a), 
there is a symplectic structure
  $\omega$ on $X$ with $[\omega] = c_f \in H^2_{\dR}(X)$, taming
  $J$ (so that $f$ is $\omega$-tame as in Definition~\ref{def:2.9}).
For fixed $J$, such structures $\omega$ form a
  convex subset of the space of all closed 2-forms.  
If only $f$ and
  $f_*J$ are fixed (otherwise allowing $J$ to vary), the space of
  such (smooth) symplectic 
structures is still contractible in any $C^k$-topology, $0
  \le k \le \infty$, or locally convex metric topology in between,
  as is the completion in such a metric. 
\item[c)] For each deformation class $\varphi$ of hyperpencils on $X$,
  there is a unique isotopy class of symplectic forms on $X$
  containing representatives $\omega$ for which some $f \in \varphi$ is 
  $\omega$-tame and $[\omega] = c_f$.
\end{itemize}
\end{thm}

\noindent
There are also versions of (a) and (b) for parametrized families of
hyperpencils; see Lemmas~\ref{lem:3.2} (as applied in proving 
Theorem~\ref{thm:2.11}(a)) and \ref{lem:3.3}, respectively.  
The completion in (b) means to complete the space of
2-forms, take the closure of the affine subspace with $[\omega] =
c_f$, and restrict to an open subset using the taming condition.
Examples include the $C^k$-spaces of taming forms, $0 \le k \le
\infty$, and many Sobolev spaces.  
The result could be stated in even
further generality; see Theorem~16 of \cite{P}.

\begin{example} \label{example:2:12} 
For a smooth algebraic variety $X \subset \complex \PP^N$, the
standard holomorphic structure on $X$ is obviously
$(\omega_{\std}|X)$-tame.  Thus the isotopy class of symplectic
structures determined by the canonical deformation class of
hyperpencils (Example \ref{example:2.6}(a)) is the one containing the
standard K\"{a}hler form $\omega_{\std}|X$.
\end{example}

\section{Proof of the Main Theorem~\ref{thm:2.11}}

Our main tool for constructing symplectic structures is the following
theorem. 
This is based on an idea, which was used by Thurston \cite{T} to construct 
symplectic structures on total spaces of surface bundles over symplectic 
manifolds, and which generalizes to bundles with symplectic fibers of 
arbitrary dimension (e.g., \cite[Theorem~6.3]{McS}). 
We now generalize still further, from bundle projections to singular maps 
suitably controlled by almost-complex structures, 
and work relative to a subset  $C$ of the domain. 
The result is general enough to apply to hyperpencils, and also 
has other applications \cite{G3}, \cite{G4}.

\begin{thm} \label{thm:3.1}
Let $f:X \to Y$ be a smooth map between manifolds, with $X -\Int C$
compact for some closed subset $C$ with a neighborhood $W_C$ in $X$.
Suppose that $\omega_Y$ is a symplectic form on $Y$, and $J$ is a
continuous, $(\omega_Y, f)$-tame almost-complex structure on $X$.  
Fix a class $c \in H^2_{\dR}(X)$.  Suppose that for each $y \in Y$,
$f^{-1}(y)$ has a neighborhood $W_y$ in $X$ containing $W_C$, with the
restriction $H^1_{\dR}(W_y) \to H^1_{\dR}(W_C)$ surjective,
and with a closed 2-form $\eta_y$ on $W_y$ such that $[\eta_y] = c|W_y
\in H^2_{\dR}(W_y)$ and such that $\eta_y$ tames $J$ on each
of the complex subspaces $\ker df_x$, $x \in W_y$. 
Suppose that these
forms $\eta_y$ all agree on $W_C$, and that the resulting form
$\eta_C$ on $W_C$ tames $J$ on $TX|C$.  
Then there is a closed 2-form
$\eta$ on $X$ agreeing with $\eta_C$ near $C$ with $[\eta] = c \in
H^2_{\dR}(X)$, and such that for all sufficiently small $t>0$
the form $\omega_t = t\eta + f^*\omega_Y$ on $X$ tames $J$ (and hence is
symplectic).
\end{thm}

To show why this theorem generalizes Thurston's construction, 
we obtain the latter as a special case. 
(See \cite{G2} for additional details.) 
Suppose $f:X \to Y$ is a fiber bundle with $X$ compact.  
Then we can take $C$ and $W_C$ to be empty.  
The required hypotheses for 
Thurston's construction (e.g., as given in \cite{McS}) include symplectic
structures on $Y$ and the fibers  $f^{-1}(y)$,  which give 
almost-complex structures on $Y$ and
the subbundle of $TX$ tangent to the fibers.  
These can easily be
combined into an almost-complex structure $J$ on $X$ that makes $f$
$J$-holomorphic and hence is $(\omega_Y, f)$-tame.  
The final hypothesis for Thurston's construction guarantees 
the existence of a suitable class $c$,  
and the remaining hypotheses of
Theorem~\ref{thm:3.1} are now easily verified (cf. \cite{G2}).  
The family of symplectic forms resulting from Theorem~\ref{thm:3.1} 
is the same one obtained by Thurston.

The proof of Theorem~3.1 is obtained by modifying Thurston's method to exploit 
the almost-complex
structure, allowing us to deal with a complicated critical set for
$f$ and work relative to a possibly nonempty subset $C$.

\begin{proof}
Fix a representative $\zeta$ of the deRham class $c$.  
For each $y \in Y$, $[\eta_y] = c|W_y$, so we can write $\eta_y = \zeta +
d\alpha_y$ for some 1-form $\alpha_y$ on $W_y$.  
Pick some $y_0\in Y$ and set $\alpha_C = \alpha_{y_0}|W_C$. 
Then for each $y$,
$d(\alpha_C - \alpha_y) = (\eta_{y_0} - \zeta) - (\eta_y -
\zeta) = 0$ on $W_C$, so $[\alpha_C - \alpha_y] \in
H^1_{\dR}(W_C)$ is defined.  
By hypothesis, any such class
extends to $H^1_{\dR}(W_y)$, so after adding a closed form on
$W_y$ to $\alpha_y$, we can assume $\alpha_C - \alpha_y$ is exact on $W_C$.  
Choosing a function $g:W_y \to \real$ with $dg =
\alpha_C - \alpha_y$ near $C$, and replacing $\alpha_y$ by $\alpha_y +
dg$, we obtain that $\alpha_y = \alpha_C$ near $C$ for each $y$. 
Since each $X-W_y$ is compact, each $y \in Y$ has a neighborhood
disjoint from $f(X-W_y)$.  
Thus, we can cover $Y$ by open sets $U_i$,
with each $f^{-1}(U_i)$ contained in some $W_y$.  Let $\{\rho_i\}$ be a
subordinate partition of unity on $Y$. 
The corresponding partition of
unity $\{\rho_i \circ f\}$ on $X$ can be used to splice the forms
$\alpha_y$; let $\eta =\zeta + d\sum_i(\rho_i \circ f)\alpha_{y_i}$. 
Clearly, $\eta$ is closed with $[\eta] = [\zeta] = c
\in H^2_{\dR}(X)$, and $\eta = \eta_C$ near $C$, so it suffices to
show that $\omega_t$ tames $J$ ($t > 0$ small). 
In preparation, perform the
differentiation to obtain $\eta = \zeta + \sum_i(\rho_i \circ f)
d \alpha_{y_i} + \sum_i (d\rho_i \circ df) \wedge
\alpha_{y_i}$. 
The last term vanishes when applied to a  pair of
vectors in $\ker df_x$, so on each $\ker df_x$ we have
$\eta = \zeta + \sum_i(\rho_i \circ f)d\alpha_{y_i} =
\sum_i(\rho_i \circ f)\eta_{y_i}$. 
By hypothesis, this is a convex
combination of taming forms, so we conclude that 
$J|\ker df_x$ is $\eta$-tame for each $x \in X$.

It remains to show that there is a $t_0 > 0$ for which
$\omega_t(v, Jv) > 0$ for every $t \in (0, t_0)$ and $v$ in the unit
sphere bundle $\Sigma \subset TX$ (for any convenient metric).  
But
$$\omega_t(v, Jv) = t \eta (v, Jv) + f^*\omega_Y(v, Jv)\ .$$ 
Since $J$
is $(\omega_Y, f)$-tame, the last term is positive for $v \notin
\ker df$ and zero otherwise.  
Since $J|\ker df$ is
$\eta$-tame, the continuous function $\eta(v,Jv)$ is positive for all
$v$ in some neighborhood $U$ of $\ker df \cap \Sigma$ in $\Sigma$. 
Similarly, for $v \in \Sigma|C$, $\eta(v, Jv) = \eta_C(v, Jv) > 0$. 
Thus, $\omega_t(v, Jv) > 0$ for all $t > 0$ when $v \in U \cup \Sigma|C$. 
On the compact set $\Sigma|(X-\Int C) - U$ containing the
rest of $\Sigma$, $\eta(v,Jv)$ is bounded and the last displayed term is
bounded below by a positive constant, so $\omega_t(v, Jv) > 0$ for $0 <
t < t_0$ sufficiently small, as required.
\end{proof}

We also need some techniques for splicing together locally defined
almost-complex structures. 
These are given by the following lemma,
whose proof appears in Section~\ref{sec:4}.  
As in Definition~\ref{def:2.2} 
(of wrapped critical points), we let $E, F \to X$ be
real vector bundles over a metrizable space, with fiber dimensions
$2n$ and $2n-2$ respectively, and this time equipped with fiber
orientations.  
We again fix a section $T:E \to F$ of $\Hom (E,F)$.  
In the applications, $T$ will be $df:TX \to f^*TY$ for some
$C^1$ map $f:X \to Y$, or a family of such maps, although the extra
generality is useful in proving the lemma below.  
We wish to splice
together locally defined complex structures on the bundle $E$, so as
to extend a preassigned structure from some closed (possibly empty)
subset $C \subset X$ to all of $X$.  
Our approach requires a
restriction on either the topology of the critical points or the
induced complex structures on the image $T(E) \subset F$.  
We assume the latter
restriction on a closed subset $D$ (which without loss of generality
contains $C$) and the former restriction elsewhere.  
The cases we
require are when $D$ equals $C$ or $X$ (allowing us to set $V$ equal to $U$ or $X$
below), but the general case poses no additional difficulties. 
A {\em 2-form\/} on $E$ or $F$ will mean a continuously varying choice of
a skew-symmetric bilinear form on each fiber of the bundle.

\begin{lem} \label{lem:3.2}
For $E,F \to X$ and $T$ as above, let $C \subset D \subset X$ be
closed subsets such that the regular points of $T|X-C$ are dense in
$X-C$, and let $\omega_F$ be a nondegenerate 2-form on $F$ (inducing
the given fiber orientation).  For some neighborhood $U$ of $C$, let
$J_C$ be an $(\omega_F, T)$-compatible complex structure on the
(oriented) bundle $E|U$.  Suppose that each $x \in X-U$ has a
neighborhood $W_x$ with an $(\omega_F, T)$-compatible complex
structure on $E|W_x$, and that for some neighborhood $V$ of $D$, these
can be chosen for all $x \in V-U$ so that the induced structures on
$T(E|W_x) \subset F|W_x$ agree with each other and with $T_*J_C$
wherever the domains overlap. Let $J_D$ denote the resulting complex
structure on the fibers of $T(E|V)$.
\begin{itemize}
\item[a)] If $n \geq 3$, assume each critical point of $T$ in $X-D$ is
  wrapped. 
Then $J_C|C$ extends to an
$(\omega_F, T)$-compatible complex structure $J$ on $E$ with $T_*J|D = J_D$. 
\item[b)] Suppose $D=X$ and $\omega_E$ is a 2-form on $E$. 
If the local complex structures on $E$ given above (including $J_C$) can be
  chosen to be $\omega_E$-tame, then we can assume $J$ is
  $\omega_E$-tame.
\item[c)] Suppose that $\partial C \subset X_0 = X - \Int C$ and (if
  $D\neq C$)
  $\partial D \subset X - \Int D$ have disjoint collar neighborhoods compatible
  with the given structures on $E$ and $F$.  
(See below.)
Then in both cases (a) and (b) above, 
the space $\J$ of all complex structures $J$ satisfying
  the given conclusions  is weakly
  contractible when $X_0$ is locally compact and contractible when
  $X_0$ is compact. 
\end{itemize}  
All of the above remains true if compatibility is replaced by taming 
everywhere, provided that $D=X$. 
\end{lem}

\noindent
Recall that $\J$ is given the compact-open topology, which equals the 
$C^0$-topology when $X_0$ is compact. 
See the proof for further details. 
The ``compatible collar'' hypothesis on $\partial C$ means that there
is a subset $K \subset X_0$ homeomorphic to $I \times \partial C$ with
$\{0\} \times \partial C$ mapping to $\partial C$ in the obvious way
and $\{1\} \times \partial C$ mapping onto $\partial K$ (the boundary in
$X_0$ in the sense of general topology), and that $E|K$ and $F|K$
can be identified with $I \times (E|\partial C)$ and $I \times
(F|\partial C)$ in such a way that $T$, $\omega_F$, $\omega_E$ (in
case (b)), and $J_D$ (if $D \neq C$) are constant over each $I \times
\{x\} \subset I \times \partial C$.

The hypothesis of wrapped critical points can be weakened, at least 
when $n=3$. 
For an open set $W\subset X$, we call an $(\omega_F,T)$-compatible complex 
structure $J$ on $E|W$ {\em $(\omega_F,T)$-extendible\/} along a 
collection of convergent sequences of regular points in $X$ if for each 
such sequence $x_i\to x$ with $x\in W$, the complex structures $T_*J$ 
on each $F_{x_i}$ (defined for all sufficiently large $i$) limit to an 
$\omega_F$-tame complex structure on $F_x$ (which is necessarily 
$\omega_F$-compatible and an extension of $T_*J$ on $T(E_x)\subset F_x$). 

\begin{adden}\label{addendum}
If $n=3$ and there are critical points in $X-D$ that are
not wrapped,  fix a sequence of regular points converging to each 
unwrapped critical point in $X-D$, and assume that the given local 
complex structures on $E$ (including $J_C)$ are $(\omega_F,T)$-extendible 
along these sequences.
Then Lemma~\ref{lem:3.2} still holds, where the structures $J$ comprising 
$\J$ are all required to be $(\omega_F,T)$-extendible along the given 
sequences, provided that in (c), the collar of $\partial C$ (if $D=C$) 
or of $\partial D$ (if $D\ne C$) contains no unwrapped critical points. 
\end{adden}

\begin{proof}[Proof of Theorem~\ref{thm:2.11}]
To prove (a),
we wish to apply Lemma~\ref{lem:3.2} to the hyperpencil $f:X-B \to
\complex\PP^{n-1}$.  First, we show that the regular points of $f$ are
dense in $X-B$ by the method previously used for Proposition~\ref{prop:2.3}. 
If any neighborhood $W \subset X-B$ consists entirely
of critical points, choose $x_0 \in W$ minimizing $\dim \ker
df_x$ and note that $\ker df$ is a smooth distribution
near $x_0$. (As before we can realize it as $\ker d(\pi \circ
f)$ for a suitable projection $\pi$.) 
By integrating a vector field in
$\ker df$, we obtain a curve of critical points in a single
fiber, contradicting finiteness in Definition~\ref{def:2.4}(3).  Now
we apply Lemma~\ref{lem:3.2} to $X-B$, with $E=T(X-B)$, $F= f^*
T\complex\PP^{n-1}$, $T=df$ and $\omega_F =
\omega_{\std}$ (pulled back to $F$).  
We let $C = D \subset X -
B$ consist of a closed, round ball (with center deleted) about each $b
\in B$ in the local charts given by Definition~\ref{def:2.4}(1), $U =
V$ be corresponding punctured open balls, and $J_C$ be the corresponding
complex structure on $E|U$.  
The required local complex structures
over neighborhoods $W_x$ exist by Definition~\ref{def:2.4}(2) at 
critical points $x$,
and are easy to construct at regular points.  
By (a) of the lemma, we
obtain the required almost-complex structure $J$ on $X-B$ (which
immediately extends over $X$), using either the given definition of a
hyperpencil or the variations of Remark~\ref{rem:2.5}(c). 
(For the last variation, set $D=X-B$.) 
Now (c) of the lemma (with the obvious
radial collar of $\partial C$) gives contractibility of  the space $\J_C$ of
$(\omega_{\std}, f)$-compatible almost-complex structures on $X$
extending $J_C$ on $C$ (with suitable $(\omega_{\std},df)$-extendability  or 
$f_*J$ fixed for the above variations).  
Let $\J = \bigcup_C \J_C$ be the space of all 
such $(\omega_{\std}, f)$-compatible
almost-complex structures on $X$ that are standard near $B$ (relative
to fixed charts as in Definition~\ref{def:2.4}(1), but on neighborhoods
of variable size).  		
Again fix $C$ and $U$ as
above (with the fixed $f_*J$ standard on $U$ in the case of the
last variation), and let $h_t: X \to X$ be a radial homotopy fixing $X-U$
with $h_0 = \id_X$ and $h_1$ collapsing $C$ into $B$.  Trivialize $TX|U$
in the obvious way, and extend $h_t$ to $\tilde{h}_t:TX \to TX$ as
$h_t \times \id_{\complex^n}$ over $U$ and $\id_{TX}$ elsewhere. 
Then pulling back induces a homotopy $H_t:\J \to \J$ with $H_0 = \id_\J$ and
$H_1(\J) \subset \J_C$.  
Composing with the previous contraction of
$\J_C$ produces the required contraction of $\J$.  
A similar argument
applies to the space of all $(\omega_{\std}, f)$-compatible
structures on $X$, since these must all agree on $TX|B$ by 
Lemma~\ref{lem:lines}(b). 
The proof of (a) is now complete.

The existence part of (b) will follow from Theorem~\ref{thm:3.1}, so
we begin by constructing suitable neighborhoods $W_y$ with forms
$\eta_y$.  Fix $C \subset U \subset X-B$ and $J$ on $X$ as above, and
a neighborhood $W_C$ of $C$ in $X-B$ with closure $\cl W_C \subset U$
a disjoint union of round balls.  
For each $y \in \complex\PP^{n-1}$, let $K
\subset F_y$ denote the (finite) subset of critical points of $f$
lying on the fiber $F_y$, and let $\Delta \subset X - U$ be a disjoint
union of closed balls, one centered at each point of $K$.  Define a
closed 2-form $\sigma$ on $\Delta \cup U$ as follows:  
Choose $\sigma$ to tame $J$ on $TX|K$.  
For $\Delta$ sufficiently small, we can then assume
$J$ is $\sigma$-tame on $\Delta$ (by openness of the taming
condition). On $U$, take $\sigma$ to be the standard symplectic form
from $\complex^n$, in the local coordinates given by (1) of 
Definition~\ref{def:2.4}, scaled so that its integral 
is $< 1/2$ on each complex
line through $0$ intersected with $U$.  
Now $J$ is $\sigma$-tame on $\Delta \cup U$.  
Since $J$ is $(\omega_{\std}, f)$-tame on $X-B$, $F_y
-K$ is a smooth (noncompact) $J$-holomorphic curve in $X - K$ whose
complex orientation agrees with its preimage orientation, and each
component intersects $B$ nontrivially (by Definition~\ref{def:2.4}(3)).  
To allow for the (presumably unlikely) possibility of $F_y$
being wildly knotted at $K$, we use the following trick: 
We can assume
$\partial \Delta$ is transverse to $F_y$, so the two intersect in a
finite collection of circles.  
Since each component of $F_y -K$
intersects $B$, we can connect each such circle to $B$ by a path in
$F_y - K$.  Let $\Delta_0 \subset \Int \Delta$ be a smaller collection
of balls surrounding $K$, disjoint from these paths and with $\partial
\Delta_0$  transverse to $F_y$. Then each component $F_i$ of the
compact surface $F_y - \Int \Delta_0$ either lies inside $\Int \Delta$
or intersects $B$. Let $W_y$ be the union of $\Int \Delta_0 \cup W_C$
with a tubular neighborhood rel boundary of $F_y - \Int \Delta_0
\subset  X - \Int \Delta_0$. Extend  each $F_i$ to a closed, oriented,
smooth surface $\hat{F}_i \subset W_y$ by arbitrarily attaching a
surface in $\Delta_0$.  Then the classes $[\hat{F}_i] \in H_2(W_y;
\integer)$ form a basis.

We now construct the required form $\eta_y$ on $W_y$ and apply 
Theorem~\ref{thm:3.1}. Since $F_y$ is $J$-holomorphic with $J$ $\sigma$-tame
on $\Delta \cup U$, $\sigma|F_i \cap (\Delta \cup W_C)$ is a positive
area form.  After rescaling $\sigma$ on $\Delta$ so that
$\int_{\hat{F}_i \cap \Delta} \sigma < 1/2$ for each $i$, we can
extend $\sigma$ over each $F_i$ intersecting $B$ as a positive area
form with $\int_{\hat{F}_i} \sigma = \#F_i \cap B$, the (positive) number of
points of $B$ in $F_i$.  
Define $\pi: W_y \to W_y$ by smoothly
splicing $\id_{W_y \cap (\Delta \cup U)} $ together with the normal bundle
projection on $W_y - \Delta$ so that $\Im \pi \subset F_y \cup \Delta
\cup U$ and $\pi|F_y \cup \Delta_0 \cup W_C$ is the identity.  Then
$\eta_y = \pi^*\sigma$ is a well-defined closed 2-form on $W_y$.
On $W_y \cap (\Delta \cup W_C)$ (away from $\partial \Delta$) $\eta_y
= \sigma$ tames $J$ (hence, each $J|\ker df_x$) and $\eta_y$ is
standard on $W_C$.  Similarly, $\eta_y|(F_y - K)$ tames $J$ on each
$T_xF_y = \ker df_x$, so after narrowing the tubular
neighborhood defining $W_y$, we can assume $\eta_y$ tames $J|\ker
df_x$ for all $x \in W_y$ (since the taming condition is open
and the set of critical points of $f$ is closed).  
For each $F_i$
intersecting $B$, we have $\langle \eta_y, \hat{F}_i \rangle =
\langle \sigma, \hat{F}_i \rangle = \#F_i \cap B = \langle c_f,
\hat{F}_i \rangle$.  
(The last equality follows, e.g.,
by computing the intersection number of $\hat{F}_i$ (pushed slightly
off $B$) with $f^{-1}(\complex\PP^{n-2}) \subset X - B$. The only
intersections occur in $U$, and each is $+1$ by the local description
of $f$ there.)  Similarly, for $F_i$ disjoint from $B$ we have
$\langle c_f, \hat{F}_i \rangle = 0 = \langle \eta_y, \hat{F}_i
\rangle$ since $\eta_y = \sigma$
is exact on the disjoint union of balls $\Delta$ containing
$\hat{F}_i$. Thus $[\eta_y] = c_f|W_y \in H^2_{\dR}(W_y)$ since
these agree on a basis of $H_2(W_y; \integer)$. 
Now we can apply Theorem~\ref{thm:3.1} 
to $X - B$ with $\omega_Y = \omega_{\std}$ on $Y =
\complex \PP^{n-1}$ and $c = c_f$. Note that $W_C$ is a disjoint union
of punctured open $2n$-balls with $n \geq 2$, so $H^1_{\dR}(W_C)
= 0$. We obtain a closed 2-form $\eta$ on $X - B$ that is standard on
$C$ (relative to the charts given in Definition~\ref{def:2.4}(1)) and
hence extends over $X$.  Then $[\eta] = c_f \in H^2_{\dR}(X)$,
and we can choose some $t > 0$ so that $\omega_t = t\eta + f^*
\omega_{\std}$ tames $J$ on $X - B$.

Unfortunately, the form $\omega_t$ is singular at $B$.  We verify this
with a local model, and find a way to eliminate the singularities.  In
the given local coordinates at $b \in B$, $J$ is the standard complex
structure on $\complex^n$, $f$ is projectivization, and $\eta$ is the
standard symplectic form. 
Up to a constant rescaling of the
coordinates, the latter can be written in ``complex spherical
coordinates'' (cf. \cite[Proposition~5.8]{McS})
as $\eta = r^2 f^*\omega_{\std} +
\frac{1}{2\pi} d(r^2) \wedge \beta$, where $r$ is the radial
coordinate on $\complex^n$ and $\beta$ is the pull-back to $\complex^n
- \{0\} \approx S^{2n-1} \times \real$ of the connection {1-form} on
$S^{2n-1}$ for the tautological bundle $S^{2n-1} \to \complex
\PP^{n-1}$, whose corresponding horizontal distribution $H$ is
orthogonal to the fibers.  (To verify this formula, note that $H$ on
$\complex^n - \{0\}$ is orthogonal to each complex line $L$ through $0
\in \complex^n$ under both the above $\eta$ and the standard
symplectic form.  
Both descend from $S^{2n-1}$ to $\omega_{\std}$
on $\CP^{n-1}$ up to a constant scale factor, and both scale by
$r^2$ radially, so they agree up to scale on $H$.  On each $L$,
$d(r^2) \wedge \beta = 2 r dr \wedge d \theta$ is standard   up to a
scale factor independent of $L$. 
The two
terms of $\eta$ are scaled compatibly since $d\eta = 0$ by
computation, using the fact that $d\beta$ pushes down  to the
tautological curvature form $2 \pi \omega_{\std}$ on
$\complex\PP^{n-1}$ with Chern class $-[\omega_{\std}] = -h$.)
Now $\omega_t(r) = (1 + tr^2)f^*\omega_{\std} +
\frac{t}{2\pi} d(r^2) \wedge \beta$ in these local coordinates.  
Clearly, this is singular at 0 since $\eta$ is nonzero there.  
However, the radial change of variables $R^2 =
\frac{1+tr^2}{1+t}$ ($t$ constant) shows that $\eta(R) = \frac{1}{1+t}
\omega_t(r)$, so there is a radial symplectic embedding
$\varphi:(\complex^n-\{0\}, \frac{1}{1+t}\omega_t) \to (\complex^n,
\eta)$ sending any deleted neighborhood of 0 to a deleted neighborhood
of the ball $R^2 \le \frac{1}{1+t}$.  
For $V \subset \complex^n$ the
image of the given coordinate chart at $b \in B$, define $\varphi_0:V
\to \complex^n$ to be a radially symmetric diffeomorphism onto an open
ball, agreeing with $\varphi$ outside of a closed ball about 0 in
$V$.  Let $\omega$ be $\varphi_0^*\eta$ near each $b \in B$ and
$\frac{1}{1+t} \omega_t$ elsewhere.  
These pieces fit together to
define a symplectic form on $X$, since $\varphi$ is a symplectic
embedding. (This construction is equivalent to blowing up $B$,
applying Theorem~\ref{thm:3.1} with $C = \emptyset$ to the resulting
singular fibration, and then blowing back down, but it avoids
technical difficulties associated with taming on the blown up base
locus.) 

To complete the existence proof for (b), we only need to verify that
the symplectic form $\omega$ on $X$ has the required properties.  Away
from $B$, we already know that $\omega = \frac{1}{1+t}\omega_t$ tames
$J$. Near $b \in B$, we have local coordinates with $J$ standard and
$\omega = \varphi_0^*\eta$, $\eta$ standard up to a constant scale
factor.  Since $\varphi_0$ is radially symmetric, it preserves the
horizontal distribution $H$ on $\complex^n-\{0\}$ and the form
$\eta|H$ up to (nonconstant) scale.  Also, $\varphi_0$ preserves each
complex line through 0, and $\varphi_0^*\eta$ is a positive area
form on each. Since these complex lines and $H$ are $\eta$-orthogonal
and $J$-holomorphic, $J$ is $\omega$-tame near $B$ and hence everywhere
on $X$.  To compute the cohomology class $[\omega] \in
H^2_{\dR}(X)$, it suffices to work outside $C$.  Then $[\omega]
= \frac{1}{1+t}[\omega_t] = \frac{1}{1+t}(t c_f +
f^*[\omega_{\std}]) = c_f$ since $[\omega_{\std}] = h
\in  H^2_{\dR}(\complex\PP^{n-1})$.  Thus, the constructed form
$\omega$ has the required properties.  

For fixed $J$, the space of taming $\omega$ with $[\omega]$ fixed is
obviously convex; the rest of the theorem depends on the following lemma.
This lemma allows us to extend a family of taming forms to a 
parametrized family of hyperpencils, starting from any subfamily with 
reasonable local topology.

\begin{lem} \label{lem:3.3}
Let $f:X - B \to Y$ be a family of hyperpencils parametrized by a
metrizable space $S$, and suppose the restriction of $f$ to some
closed subset $S_0 \subset S$ is $\omega_0$-tame for some continuous
family $\omega_0$ of continuous  2-forms $\omega_s$ on the fibers
$X_s$, $s \in S_0$.  
Suppose $S_0$ has a neighborhood $U_0 \subset S$
with a $C^0$ retraction $r: U_0 \to S_0$ covered by retractions $r_X:
\pi_X^{-1}(U_0) \to \pi_X^{-1}(S_0)$ and $r_Y:\pi_Y^{-1}(U_0) \to
\pi_Y^{-1}(S_0)$ with fiberwise derivatives $dr_X$ and $dr_Y$
continuous over $U_0$, $r_Y$ preserving $\omega_{\std}$, $r_X$ a
fiberwise diffeomorphism preserving $B$ 
and complex linear near $B$ with respect to the vector bundle
structure, and $f \circ r_X = r_Y \circ f$.  Then $\omega_0$ extends
to a family $\omega$ parametrized by $S$, taming $f$.  
If $\omega_0$ consists of $C^k$-forms  varying continuously in the 
$C^k$-topology, 
$(0 \le k \le \infty)$, then $\omega$
inherits this property (provided $dr_X$ is $C^k$-continuous), and if
the family $\omega_0$ is $C^k$ 
(for $\pi_X:X\to S$ a $C^{k+1}$-bundle map of  $C^{k+1}$-manifolds) then so is $\omega$.
Similarly, closure of each $\omega_s$ and the condition $[\omega_s] =
c_f \in H^2_{\dR}(X_s)$ are inherited by $\omega$.  
We can also assume
that $\omega$ tames a preassigned $J$ on $T^v X$ as specified in 
Lemma~\ref{lem:2.10} (provided that $\omega_0$ tames $J$ over $S_0$),  
or that $\omega$ tames some $J$ for which $f_* J$ on $\Im df\subset 
f^* T^vY$ is preassigned, $r_Y$ invariant, standard near $B$, 
$\omega_{\std}$-compatible,  
and induced from local complex structures on $T^vX$ as in 
Definition~2.9 (over $S_0$) and Definition~2.7(2) (over $S$). 
\end{lem}

\begin{proof}
If $J$ was not preassigned, construct it as follows: 
Find an $\omega_0$-tame $J_0$ over
$S_0$ by Lemma~\ref{lem:2.10}, and pull it back over $U_0$ by $r_X$.
Then $J_0$ is continuous and $(\omega_{\std}, f)$-compatible
over $U_0$, standard near $B$ there, and if $f_* J$ is preassigned,
it agrees with $f_* J_0$ there. 
Now $J_0|\pi_X^{-1}(S_0)$ extends to $J$ on $T^v X$ as specified in  
Lemma~\ref{lem:2.10} (where $\omega_0$-tameness of $J$ only
applies over $S_0$), by Lemma~\ref{lem:3.2}(a) applied to $X - B$
with $E = T^v(X-B)$, $F = f^* T^v Y$, $C = \pi_X^{-1}(S_0)$ union
the closure of a suitable neighborhood of $B$, and $D$ equal to $C$ or $X$,
depending on whether $f_* J$ was given.  
(Note that $X$ is metrizable since it is locally metrizable and
paracompact Hausdorff, 
e.g., \cite[\S42]{Mu}, \cite{P}.)
By the existence part of
(b) of Theorem~2.11, each $X_s$, $s \in S - S_0$, has a smooth
symplectic form $\omega_s$ taming $J|X_s$, with $[\omega_s] = c_f$. 
By openness of the taming condition and compactness of each fiber $X_s$,
each $\omega_s$ extends over some neighborhood $W_s$ of $s \in S -
S_0$ (e.g., via local triviality of $\pi_X$) so that it tames
$J|\pi_X^{-1}(W_s)$. Similarly, each $s_0 \in S_0$ has a neighborhood
$W_{s_0}$ in $S$ over which $\omega_0$ extends as a taming form,
preserving any additional conditions. 
(If only a continuous family is required, pull back by $r_X$.  
For a $C^k$-family of closed forms,
$k\ge 1$, locally trivialize, $C^k$-extend the map
$[\omega_s]:S_0 \to H^2_{\dR}(X_s)$ over some $W_{s_0}$
(possibly with $[\omega_s] = c_f$ by hypothesis), find a
$C^k$-family of representatives, and correct by $d\alpha$ for a
suitably extended family of 1-forms $\alpha$ to recover the original
subfamily $\omega_0$.) 
Using a partition of
unity on $S$ subordinate to the cover $\{W_s|s \in S\}$, splice
together the local families $\omega_s$ into a global family $\omega$
extending $\omega_0$.  
The taming and closure conditions are preserved
since each $\omega|X_s$ is a convex combination (with constant
coefficients) of such forms.  
The lemma follows immediately.
\end{proof}

To complete the proof of Theorem~2.11, first let $\omega_s$, $s =
0,1$, be symplectic forms on $X$ associated to hyperpencils $f_s$ in a
deformation class $\varphi$ as in (c). 
Then there is a deformation agreeing with $f_0$ and $f_1$,
respectively, on neighborhoods of $0,1$ in $I$.  
The above lemma, with
$S=I$ and $S_0 = \{0,1\}$, gives a smooth family 
of symplectic forms $\omega_s$, $0 \leq s
\leq 1$, interpolating between $\omega_0$ and $\omega_1$, with
$[\omega_s] = c_{f_0} = c_{f_1}$ for each $s$.  
By Moser's Theorem
\cite{M}, any deformation of cohomologous symplectic forms is realized
by an isotopy, proving (c).  
A similar argument proves the last
sentence of (b), completing the proof of the theorem: 
For fixed $k$, $1 \le k\le \infty$, let $\C$ denote the $C^k$-space of 
all $C^k$-symplectic forms on $X$ as in (b) (for fixed $f$ and $f_* J$).  
Any continuous map
$\varphi: \partial D^{m+1} \to \C$ ($D^{m+1}$ an $(m+1)$-disk) 
can be interpreted as a $C^k$-continuous family $\omega_0$ of symplectic 
forms on $X$ parametrized by $\partial D^{m+1}$.  
Applying the above lemma
to the constant family of hyperpencils on $X$ with $S=D^{m+1}$, $S_0 =
\partial D^{m+1}$, and the given $f_* J$ (which ensures that $f$ over
$S_0$ 
is $\omega_0$-tame), we extend to a family
parametrized by $D^{m+1}$, or equivalently a continuous map $D^{m+1}\to \C$
extending $\varphi$.  
This shows that $\C$ is weakly contractible ($\pi_m(\C) = 0$ for all $m$).  
But $\C$ is an open subset of the affine
subspace of 2-forms determining $c_f \in H^2_{\dR}(X)$, so 
it is a metrizable manifold of infinite dimension.  
In particular, it is an ANR, so weak contractibility
implies contractibility \cite{P}.  
Since the $C^\infty$-forms are dense in the $C^0$-space of forms, 
and $\C$ extends to an open subset of a closed affine subspace of the 
latter, the required assertion follows from Theorem~16 of \cite{P}. 
\end{proof}

\section{Proofs of Lemmas  \ref{lem:2.10} and \ref{lem:3.2}} \label{sec:4}

The proofs of Lemmas~\ref{lem:2.10} and \ref{lem:3.2} require a 
canonical method for interpolating between almost-complex structures.  
Our approach is a generalization of that of \cite{ABKLR}, Proposition~6.2.  
This depends on the $r=-\frac{1}{2}$ case of the following 
proposition, a similar form of which is stated without proof on p.100 
of  \cite{ABKLR}.  
The proof below follows a suggestion of L.~Sadun.

\begin{prop} \label{prop:4.1}
For $r \in \real$, let $\rho_0: \complex - (-\infty,0] \to \complex$
be the branch of $z^r$ with $\rho_0(1)=1$.  
Let $\A \subset \GL(m, \complex)$ 
be the open subset of matrices  with no eigenvalues in $(-\infty,0]$.
Then there is a unique holomorphic map $\rho: \A\to \GL(m, \complex)$, 
which will be denoted by $\rho(A)=A^r$, such that each $\lambda$-eigenvector
of $A$ is a $\rho_0(\lambda)$-eigenvector of $A^r$.  
This map has the following properties:
\begin{itemize}
\item[a)] For $n \in \integer, A^n$ agrees with its usual meaning.
\item[b)] For $|r| \leq 1$ or $s \in \integer, (A^r)^s = A^{rs}$.
\item[c)] If $T:\complex^m \to \complex^k$ is a linear transformation with
  $TA=BT$, then $TA^r = B^rT$ (whenever both sides are defined).
\item[d)] If $A$ is real, then so is $A^r$.
\item[e)] If $A$ is real and self-adjoint with respect to a given inner
  product $g$ on $\real^m$, then so is $A^r$.
\item[f)] For $n \in \integer -\{0\}$, $A^{1/n}$ is the unique solution
  to the equation $X^n = A$ for which all eigenvalues of $X$ lie in 
  $\Im\rho_0$.
\end{itemize}
\end{prop}

\noindent
Note that (d) and (c) imply that $\rho$ is canonically defined for 
any finite-dimensional real vector space. 

\begin{proof} 
For any $A \in \A$, the Jordan form of $A$ splits $\complex^m$ as the 
direct sum of the generalized eigenspaces $V_\lambda = 
\ker (A - \lambda I)^m$ ($\lambda$ ranging over the eigenvalues of $A$). 
On each $V_\lambda$, $A$ has the form $\lambda(I+N)$ for
the nilpotent transformation $N = \frac{1}{\lambda}A-I$. 
Set $\rho(A) =
\rho_0(\lambda) p(N)$ on $V_\lambda$, where $p(z)$ is the power
series expansion about 0 of the function $\rho_0(1+z)=(1+z)^r$.
(Note that $p(N)$ is a polynomial in the nilpotent transformation $N$.) 
On each $\lambda$~-~eigenspace, $N=0$ and $A^r=\rho_0(\lambda) I$
as required.  
Properties (a) and (b) follow immediately from the 
corresponding properties for $\rho_0$ and $p$.  
Property (c) is also 
immediate, once we observe that the condition $TA=BT$ implies $T$
preserves or annihilates each generalized $\lambda$-eigenspace. 
For (d), take $A$ real and note that $\real^m$ is the span of the 
vectors $v + \bar{v}, v \in V_\lambda$ (where $\lambda$ ranges 
over all eigenvalues of $A$).  
Since $\bar{v} \in V_{\bar{\lambda}}$, 
$A^r(v+\bar{v}) = A^rv + A^r\bar{v} = A^rv+\overline{A^rv}$, so this is 
real as required.  
Now (e) follows immediately from $g$-orthogonal diagonalizability. 
For (f), take any $X$ as given and write it as
$\lambda(I+N)$ on each of its generalized eigenspaces $V_\lambda$. 
Then $A = X^n = \lambda^n(I+N)^n$ on each $V_\lambda$, and the last
factor has the form $I+(\text{nilpotent})$.  
Thus $A^{1/n} =
\rho(X^n) = X$, since $\lambda$ is the unique $n^{th}$ root of 
$\lambda^n$ in $\Im\rho_0$ and the power series for $\rho_0(1+z)$ 
inverts the exponentiation of $(I+N)^n$.

It remains to prove holomorphicity of $\rho$, from which 
uniqueness follows by density of diagonalizable matrices in $\A$.
First note that the equation $\det(A-\lambda I) = 0$ in 
$(A, \lambda)$ exhibits the eigenvalues of matrices in $\A$ as an 
algebraic variety in $\A\times\complex$.  
The subset $\cS \subset \A$
of matrices failing to have $m$ distinct eigenvalues is then also a
variety, and over $\A - \cS$ the $m$ eigenvalues vary holomorphically.
One can now locally construct holomorphically varying bases of 
eigenvectors over $\A -\cS$, and in these bases, $\rho$ is easily
seen to be holomorphic as required.  
To show that $\rho$ extends from 
$\A - \cS$ to some holomorphic map $\hat{\rho}:\A \to \text{GL}(m,\complex)$, 
note that nondiagonalizable matrices with $m-1$ distinct eigenvalues 
form a dense open subset of $\cS$. 
Given such a matrix $A$, restrict it to its 
2-dimensional generalized eigenspace $V_\lambda$, and note that in a 
suitable basis we obtain the matrix $A_0$ in the family
\[
A_z = 
\begin{pmatrix} 
\lambda+z & 0 \\[6pt]
1         & \lambda 
\end{pmatrix}, 
\text{ with } 
\rho(A_z) = 
\begin{pmatrix} 
\rho_0(\lambda+z) & & 0 \\[6pt]
\dfrac{\rho_0(\lambda+z)- \rho_0(\lambda)}{z} & & \rho_0(\lambda) 
\end{pmatrix}
\]
for small $z \neq 0$. 
(The latter equality is easy to verify using the
eigenvectors $(0,1)$ and $(z,1)$ for $A_z$.)  
Clearly, $\rho(A_z)$ extends holomorphically over $z=0$. 
The Cauchy Integral Formula now
provides the required holomorphic extension $\hat{\rho}$ on $\A$
(since the remaining subset of $\cS$ has higher 
codimension in $\A$; cf.\ Hartogs' Theorem \cite{GH}).

Finally, we show $\hat{\rho} = \rho$.  
First, we assume $r = 1/n, n \in \integer - \{0\}$. 
For $A \in \cS$, consider a sequence $(A_i)$ in $\A - \cS$ converging to it.  
Then $\rho(A_i) = \hat{\rho}(A_i) \to \hat{\rho}(A)$. 
By (a) and (b), $(\rho(A_i))^n = A_i$; taking the limit shows that
$(\hat{\rho}(A))^n = A$. 
Similarly, this sequence allows us to write each eigenvalue $\lambda$
of $\hat{\rho}(A)$ as $\lim\lambda_i$ for $\lambda_i$ an 
eigenvalue of $\rho(A_i)$.  
By the definition of $\rho$, $\lambda_i\in \Im \rho_0$, so 
$\lambda_i = \rho_0(\lambda_i^n)$ and $\lambda=\rho_0(\lambda^n)$.
Now $\hat{\rho}(A) = \rho(A)$ by (f) as required.  
If $r = p/q \in \que$, consider $A_i \to A$ as before. 
Then $\hat{\rho}(A_i) = (A_i^{1/q})^p$, so
$\hat{\rho}(A) = (A^{1/q})^p = \rho(A)$  (by the previous case).  
The case of irrational $r$ now follows by 
continuity of $\rho$ and $\hat{\rho}$ with respect to $r$.
\end{proof}

\begin{cor} \label{cor:4.2}
For a real, finite-dimensional vector space $V$, let
$\B \subset \Aut(V)$ be the open set of linear operators
with no real eigenvalues, and let $\J \subset \B$ denote 
the set of complex structures on $V$  (for both orientations of $V$).  
Then there is a canonical real-analytic retraction $j:\B \to \J$.  
For any linear transformation $T:V \to W$ with 
$TA = BT$, we have $Tj(A) = j(B)T$ (whenever both sides are defined).
\end{cor}

\begin{proof}
We generalize \cite{ABKLR}.  For $B \in \B$, note that $-B^2$ 
has no real eigenvalues $\leq 0$. 
(If $-\lambda^2$ were an 
eigenvalue of $-B^2$ with $\lambda \in \real$ then $\pm \lambda$
would be an eigenvalue of $B$, since $0 = \det(-B^2 + 
\lambda^2 I) = \pm \det(B + \lambda I)\det(B - \lambda I)$.)
Thus, we can define $j(B)$ to be $B(-B^2)^{-1/2}$.  
The two factors commute by Proposition~\ref{prop:4.1}(c), so
$(j(B))^2=-I$ as required.  
For $J \in \J$, $j(J)=J$.  
The rest of the corollary also follows immediately.
\end{proof}

Now fix $J_1, \dots, J_k \in \J$ and let $t=(t_1, \dots, t_k)$ 
vary over the simplex $\sum t_i = 1$, each $t_i \geq 0$. 
Suppose each $B_t = \sum t_iJ_i$ has no real eigenvalues.  
For example,
this is guaranteed if $J_1, \dots, J_k$  are all $\omega$-tame for a 
fixed $\omega$ on $V$ (since for $v \neq 0$, $\omega(v, B_t v) =
\sum t_i \omega(v, J_i v) > 0$ but $ \omega(v,v) = 0$). 
Then we obtain an analytic simplex $j(B_t)$ of complex structures
(complex simplex?) with vertices  $J_1, \dots, J_k$. 
We show that if the vertices are $\omega$-compatible for a fixed $\omega$ on
$V$, then so is each $j(B_t)$ (cf. \cite{ABKLR}): Compatibility
implies that the bilinear form $g_i(v,w) = \omega(v, J_iw)$ is
positive definite and symmetric for each $i$, as is the form 
$g = \sum t_i g_i$. 
Since $g(v,w) = \omega(v, B_t w)$, we have
$g(B_t v, w) = \omega(B_t v, B_t w) = - \omega(B_tw, B_t v) = 
-g (B_tw, v) = - g(v, B_t w)$.  
Thus $B_t$ is skew-adjoint with
respect to the inner product $g$, so $-B_t^2$ and 
$(-B_t^2)^{-1/2}$ are self-adjoint (the latter by 
Proposition~\ref{prop:4.1}(e)).  
Hence $j(B_t)$ is skew-adjoint, 
implying that $\omega(j(B_t)v, j(B_t)w) = g(j(B_t)v, B_t^{-1}j(B_t)w)
= g(v, B_t^{-1}w) = \omega(v,w)$ as required.  
Furthermore, $\omega(v, j(B_t)v) = \omega(v, B_t(-B_t^2)^{-1/2}v) =
g(v, (-B_t^2)^{-1/2}v) >0$ for $v \neq 0$, since $-B_t^2$ 
and hence $(-B_t^2)^{-1/2}$ are $g$-orthogonally 
diagonalizable with all eigenvalues positive.

\begin{ques} \label{ques:4.3}
If $J_1, \dots, J_k$ are only given to be $\omega$-tame 
(for a fixed $\omega$), is each $j(B_t)$ $\omega$-tame?
\end{ques}

\noindent 
An affirmative answer would allow us to replace compatibility 
by taming throughout the paper, 
define $\omega$-tame
hyperpencils in a purely local way, remove $f_*J$ from the statement 
of the Main Theorem~2.11(b), and significantly simplify
some of the proofs.  
For example, (b) of Lemma~\ref{lem:3.2} with
$D=C$ would follow immediately from the proof of (a), rendering 
the case of the lemma with $D\ne C$ unnecessary, along with the nonlocal 
condition on $f_*J_\alpha$ in the definition of $\omega$-tame hyperpencils, 
and the fixed $f_*J$ in Theorem~2.11(b).  

\begin{proof}[Proof of Lemma~\ref{lem:2.10}]
The main difficulty is that Lemma~\ref{lem:3.2}(b) does not apply
directly at the base locus.  
For an $\omega$-tame hyperpencil
$f:X-B \to \complex\PP^{n-1}$ and $b \in B$, Definition \ref{def:2.9} gives a 
neighborhood $W_b$ of $b$ in $X$ with an $\omega$-tame $J_b$
on $W_b$ for which $J_b|W_b - \{b\}$ is 
$(\omega_{\std}, f)$-compatible.  
If $f_*J_b$ is not standard near $b$, we must correct this by a perturbation. 
First note that $J_b|T_bX$ must agree with the standard structure $i$, 
by (b) of the following lemma, which is proved below.

\begin{lem}\label{lem:lines}
{\rm (a)} A linear complex structure $J$ on $\real^{2n}$, $n\ne 1$, 
is determined by its 1-dimensional (oriented) complex subspaces. 

{\rm (b)} If $f: \complex^n-\{0\}\to \complex \PP^{n-1}$ denotes 
projectivization, $n\ge 2$, and $J$ is a continuous (positively
oriented) almost-complex 
structure on a neighborhood $W$ of $0$ in $\complex^n$, with $J|W-\{0\}$ 
$(\omega_{\std},f)$-tame, then $J|T_0\complex^n$ is the standard complex 
structure.
\end{lem}

\noindent
Thus there is a neighborhood $W$ of $b$ in $W_b$ for which the operators
$A_t=(1-t)J_b+ti$, $0 \leq t \leq 1$, have no real 
eigenvalues, and for which the resulting complex structures 
$j(A_t)$  from Corollary~\ref{cor:4.2} are $\omega$-tame.  
Since $f_* j(A_t) = j((1-t)f_* J_b + tf_* i)$ by the naturality 
statement in Corollary~\ref{cor:4.2}, the paragraph 
following the proof of that corollary shows that each $j(A_t)$
on $W-\{b\}$ is $(\omega_{\std}, f)$-compatible.  
Let $\rho:X \to [0,1]$ be a continuous function with support in 
$W$ and identically 1 near $b.$  
Then we can add 
$(W, j(A_\rho))$ to the open cover $\{W_\alpha\}$ and replace
each $W_\alpha$ by $W_\alpha-\supp \rho$, 
recovering the hypotheses of the lemma, with the unique
structure near $b$ standard.  
After applying this 
procedure when necessary to  each $b \in B$, we invoke Lemma~\ref{lem:3.2}(b),
with $E = T(X-B)$, $F = f^* T \complex \PP^{n-1}$, $T = df$, and $C$ 
equal to the closure of a sufficiently small neighborhood of
$B$ so that we can let $J_C$ be the standard structure $i$
near each $b \in B$.  
The resulting structure on $X-B$ extends 
in the obvious way over $B$, completing the proof for a single hyperpencil. 

A similar argument applies to families over metrizable parameter spaces.  
Cover $B \subset X$ by neighborhoods $W$ as before, and
choose $\rho:X \to [0,1]$ supported in their union.  
The resulting forms $f_* j(A_\rho)$ fit together as required, and 
Lemma~\ref{lem:3.2}(b) again completes the proof. 
\end{proof}

\begin{proof}[Proof of Lemma~\ref{lem:lines}]
To reduce (b) to (a), note that  
each complex line through $0$ (with
respect to the standard complex structure $i$ on $\complex^n$)
intersects $W$ in  
a $J$-complex curve (with the same orientation), since
its tangent spaces away from $0$ are given by $\ker df$ and 
$J$ is $(\omega_{\std}, f)$-tame.  
($J$-complexity extends over $0$ by continuity.)  
Thus, the $i$-complex lines through
$0$ in $T_0\complex^n$ are also $J$-complex lines.  
To prove (a), assume each 1-dimensional $i$-linear subspace of 
$\complex^n =\real^{2n}$ is $J$-linear,  
pick two $i$-linearly independent vectors 
$v,w \in \complex^n$ and let $W = \Span_iw$ and $V_\lambda =
\Span_i(v+\lambda w)$ for $\lambda \in \complex$.  
These subspaces are both $i$-complex and $J$-complex lines. 
Thus, the projections of $V_\lambda \subset V_0 \oplus W$  to
$V_0$ and $W$ determine a map $\varphi_\lambda: V_0 \to W$ (whose 
graph is $V_\lambda$) that is both $i$-linear and $J$-linear,
and is determined by the condition $\varphi_\lambda (v) = \lambda w$.
Now $\psi = \varphi_1^{-1} \circ \varphi_i:V_0 \to V_0$ is both $i$-
and $J$-linear and given by $\psi(v)=iv$, so $\psi$ is  
multiplication by $i$.  
It follows that $J$ commutes with
$i$ on the complex line $V_0$, so that these two complex 
structures on $V_0$ agree.  
But $V_0$ was chosen arbitrarily, so $J=i$ everywhere.
\end{proof}

\begin{proof}[Proof of Lemma~\ref{lem:3.2} and Addendum~\ref{addendum}]
We begin proving (a) with a partition of unity argument 
that works outside of $D-C$. 
Cover the space $X$ by open sets $W_\alpha$ with complex
structures $J_\alpha$ on $E|W_\alpha$ as in the lemma. 
We can  assume that $U$ is the unique  $W_\alpha$ intersecting
$C$ (with corresponding $J_\alpha = J_C$), and that 
any $W_\alpha$ intersecting $D$  
lies in $V$, so that $T_* J_\alpha = J_D|W_\alpha$ there. 
Since $X$ is metrizable, it is 
paracompact by Stone's Theorem (e.g., \cite{MS}, \cite[Theorem~41.4]{Mu}),
so there is a partition of unity $\{\rho_\alpha\}$ 
subordinate to the covering $\{W_\alpha\}$. 
Let $A = \sum \rho_\alpha J_\alpha: E \to E$.  
To make this into
a complex structure by the retraction in Corollary~\ref{cor:4.2}, 
we must investigate where $A$ could have real eigenvalues. 
First note that the map $B = \sum \rho_\alpha T_* 
J_\alpha : T(E) \to T(E)$ has no real eigenvalues on 
any fiber since each $T_* J_\alpha$ is $\omega_F$-tame.  
Since $TA=BT$, each 
$\lambda$-eigenvector of $A$ is mapped by $T$ to $0$ or
a $\lambda$-eigenvector of $B$, so any real eigenvector
of $A$ must lie in $\ker T$.  

To rule out real eigenvectors in $\ker T$, 
recall the subset $P \subset \ker T \subset E$ introduced  for 
Definition~\ref{def:2.2} (of wrapped points). 
For each $x \in X$, any $v \in P_x$ can be written as $\lim v_i$
with $v_i \in \ker T_{x_i}$ for some sequence 
$(x_i)$ of regular points converging to $x$.  
After passing to a subsequence, we can assume the oriented 2-planes
$\ker T_{x_i}$ (in the preimage orientation) converge to an 
oriented 2-plane $\Pi \subset P_x$ containing $v$. 
For each $J_\alpha$ defined on $E_x$, each $\ker T_{x_i}$  ($i$ large) 
will be $J_\alpha$-complex (compatibly oriented) as will
the limit $\Pi$.  
It is now easy to construct a decomposition 
$\Span_\real P_x = \bigoplus \Pi_j$, where each 
oriented real 2-plane $\Pi_j$ is a $J_\alpha$-complex line 
for each $J_\alpha$ defined at $x$.  
Clearly, the quotient $Q_x = \ker T_x / \Span_\real P_x$ inherits
a complex structure $\bar{J_\alpha}$ from each such 
$J_\alpha$, and these are all compatible with the same
orientation on $Q_x$ (inherited via $\ker T_x$ from
$\omega_F | T(E_x))$.  
If $x$ is wrapped, then
$\dim_\complex Q_x \leq 1$, so each $\bar{J_\alpha}$ is 
$\omega$-tame for some fixed $\omega$ on $Q_x$. 
It follows as 
before that $\sum\rho_\alpha(x)\bar{J_\alpha}$ has no real
eigenvalues on $Q_x$, so any real eigenvector of $A_x$ 
lies in $\Span_\real P_x = \bigoplus \Pi_j$.  
But a direct sum $\omega$ on this space tames each $J_\alpha$, so we 
conclude that $A_x$ has no real eigenvalues when $x$ is wrapped.  
If $x \in X-D$ is not wrapped (Addendum~\ref{addendum}),  
then by hypothesis $n=3$, 
and $\dim_\complex Q_x \geq 2$, so $T_x = 0$ (since 
density of regular points implies $\dim_\complex \Span_{\real}P_x \geq 1$).  
We are also given a sequence $x_i\to x$ of regular points along which 
each $J_\alpha$ is $(\omega_F,T)$-extendible, so each relevant $T_*J_\alpha$ 
extends continuously and $\omega_F$-compatibly to $F|\{x_i\} \cup\{x\}$. 
If we pass to a suitable subsequence, $\frac{T_{x_i}}{\|T_{x_i}\|}$ will 
converge to a nonzero transformation that is complex linear for each 
relevant $\alpha$. 
Then over the subspace $\{x_i\} \cup \{x\}$, $\frac{T}{\|T\|}$ will 
have a wrapped point at $x$, and each relevant $J_\alpha$ will be 
$(\omega_F,\frac{T}{\|T\|})$-compatible, so 
the previous argument again shows that $A_x$ has no real eigenvalues. 
Thus, $A$ can have real eigenvalues only  
at unwrapped critical points in $D$.  

Corollary~\ref{cor:4.2} provides a complex structure $j(A)$ 
that  is well-defined 
and continuous except where $A$ has real eigenvalues. 
By construction, $j(A) = J_C$ on
some neighborhood of $C$, and $T_* j(A) = j(B)$, which 
is defined everywhere and equals $J_D$ on a neighborhood $V'$ of $D$.  
Furthermore $j(A)$ is $(\omega_F, T)$-compatible (in particular, 
$(\omega_F, T)$-tame) since each $T_* J_\alpha$, and 
hence $j(B)$, is $\omega_F$-compatible.  
(This is the
only place where we require compatibility instead of
taming; cf. Question \ref{ques:4.3}. 
Our use of compatibility in proving Lemma~2.10 can be avoided using 
openness of taming, cf. \cite[second paragraph following Addendum~2.6]{G3}.) 
Similarly, in the case of Addendum~\ref{addendum}, 
$j(A)|X-D$ is $(\omega_F,T)$-extendible along the given 
sequences,  as is any structure mapping to $j(B)$. 
This proves (a) of Lemma~\ref{lem:3.2} (and its addendum) when $D=C$.
The general case reduces to the case $D=X$ of the lemma after we 
intersect each $W_\alpha$ with $V'$, add a new 
$W_\alpha = X - D$ with $J_\alpha = j(A)|X-D$, and
extend $J_D$ to $j(B):T(E) \to T(E)$.

We now complete the proof of  (a) in the remaining case
$D=X$, while simultaneously preparing for (b). 
Recall that we just showed each $v \in P$ lies
in an oriented 2-plane $\Pi = \lim \ker T_{x_i} 
\subset P_x$, for some sequence $(x_i)$ of regular 
points converging to $x$, and that $\Pi$ is a complex
line for each $J_\alpha$ (including $J_C$) defined on $E_x$.  
Each such $J_\alpha$ now induces a complex structure $J_\Pi$
on the quotient $E_x/\Pi$, and this structure is the 
limit of the corresponding structures on $E_{x_i}/\ker T_{x_i}$. 
Since these latter structures are also 
determined by $J_D$ on $T(E_{x_i})$ via the 
isomorphism induced by $T_{x_i}$, they are 
independent of $\alpha$, as is $J_\Pi$.  
If $Z \subset X$ denotes the subset of $X$ for which 
$\dim_\complex\Span_{\real} P_x \geq 2$, then 
each $E_x$ with $x \in Z$ contains at least two
such distinct planes $\Pi_1$ and $\Pi_2$, and the
$J_\alpha$-complex monomorphism $E_x \hookrightarrow
E_x/\Pi_1 \oplus E_x/\Pi_2$ shows that the structures
$J_\alpha$ all agree on $E_x$.  
Thus, over $C' = C 
\cup \cl Z$, the structures $J_\alpha|W_\alpha
- C$ and $J_C$ all fit together into a continuous
$(\omega_F, T)$-compatible structure $J_{C'}$ on $E|C'$ with
$T_* J_{C'} = J_D|C'$, depending only on $J_C$ and $J_D$. 
We can now construct a $J_{C'}$-Hermitian
fiber metric on $E|C'$ whose real part extends
to a fiber metric $g$ on $E$ (using partitions of unity
and the Tietze Extension Theorem). 
Over $X-C'$, $P$ is an oriented (real) 2-dimensional subbundle of $E$. 
(To verify that $P|X-C'$ is a continuous section
of the Grassmann bundle, note that for $x \in X-C'$
and any neighborhood $V_\Pi$ of the plane $\Pi = P_x$,
there is a neighborhood $U_x$ of $x$ in $X$ on which
all regular points map into $V_\Pi$; otherwise 
one could construct a sequence $x_i \to x$ with
$\ker T_{x_i}$ converging to a plane $\neq \Pi$.
Now $U_x - C'$ maps into the closure of $V_\Pi$.) 
Using the given orientation of $P|X-C'$ and the metric $g$,
define $J|X-C'$ to be counterclockwise 
$\frac{\pi}{2}$-rotation on $P$ and the structure 
determined by $J_{P_x}$ on each $P_x^\perp \cong
E_x/P_x$. (Note that  $J_{P_x}$ is continuous in 
$x$ since it is locally induced by the structures $J_\alpha$.) 
For each $x \in X-C'$, $T_x$ factors 
through $E_x/P_x$ on which $J$ and $J_\alpha$ agree, 
so $J|X-C'$ is $(\omega_F, T)$-compatible, with 
$T_* J|X-C' = J_D|X-C'$.  
Thus the proof of (a) is completed by showing that $J|X-C'$ and 
$J_{C'}$ fit together continuously at each $x \in C'$.  
If this fails, then there is  a sequence $x_i \to x$ 
for which each $J_{x_i}$  lies outside a fixed 
neighborhood of $J_x$ in the bundle $\End (E)$. 
By continuity of $J_{C'}$, we can pass to a subsequence 
$(x_i)$ lying in $X-C'$, and then further assume the 
oriented 2-planes $P_{x_i}$ converge to some $\Pi$ 
at $x$ by compactness of the Grassmann manifold. 
Since $g$ is continuous and Hermitian at $x$, 
it is now routine to obtain the contradiction that $J_{x_i} \to J_x$.
(Compare with a fixed $J_\alpha$ on the given orthogonal summands.) 

To prove (b) of the lemma, we introduce a 2-form
$\omega_E$ on $E$ that is assumed to tame each 
$J_\alpha$ (including $J_C$).  
On $X-C'$, let $Q$ denote the $\omega_E$-orthogonal complement 
of the 2-dimensional (real) oriented subbundle $P \subset E|X-C'$. 
Since each $P_x$ is a complex line for some $\omega_E$-tame $J_\alpha$ 
on $E_x$, $\omega_E$ is nondegenerate (and positively 
oriented) on $P$.  
Thus $E|X-C' = Q \oplus P$ is an $\omega_E$-orthogonal direct sum splitting.  
Now any subbundle of $E|X-C'$ complementary to $P$ can be 
written as $\graph \psi$ for some continuous section $\psi$ of $\Hom (Q,P)$. 
For $J$ as constructed in the 
previous paragraph, let $J_\psi$ denote the complex 
structure on $E|X-C'$ given by $J$ on $P$ and by
$J_{P_x}$ on each $\graph \psi_x \cong E_x/P_x$.  
To express any $J_\psi$ in terms of $J_0$ (which is given
by $J_{P_x}$ on each $Q_x$), note that these agree on the 
quotient $E_x/P_x$, so for any $(q, \psi(q)) \in 
\graph \psi \subset Q \oplus P$ we have $J_\psi(q, \psi(q))
= (J_0 q, \psi(J_0 q))$.  
Then for any $(q,p) \in Q 
\oplus P$ we have $J_\psi(q,p) = {J_\psi(q, \psi(q)) +
J_\psi(0, -\psi(q) + p)} = (J_0 q, (\psi J_0 - J\psi)q + Jp)$.  
In particular,
\begin{equation} \label{eq:4}
\omega_E((q,p), J_\psi(q,p)) = \omega_E(q, J_0q) +
\omega_E(p,Jp) + \omega_E(p, (\psi J_0 - J\psi )q).
\end{equation}
Note that $J$ is $\omega_E$-tame on $P|X-C'$ since 
the fibers are correctly oriented $J$-complex lines.  
In the next paragraph, we will show that $J_0$ is 
$\omega_E$-tame on $Q$, so the first two terms on the
right side of Equation~\eqref{eq:4} have positive sum 
whenever $(q,p) \neq (0,0)$.  
Now choose $\psi$ so that
$\graph \psi$ is a $J$-complex subbundle of $E|X-C'$.
Then $J_\psi = J|X-C'$. 
At each $x \in C'$, $J$ 
agrees with some $J_\alpha$, so it is $\omega_E$-tame on $E_x$.  
By openness of the taming condition, we 
conclude that $J$ is $\omega_E$-tame over some 
neighborhood $U'$ of $C'$.  
Thus, the left side of 
\eqref{eq:4} is positive for $x \in U'-C', (q,p) \ne (0,0)$.  
Replacing $\psi$ by $\rho \psi$ in \eqref{eq:4}
for any $\rho:U'-C' \to [0,1]$, we obtain a convex
combination of two positive quantities, showing that
$J_{\rho \psi}$ is $\omega_E$-tame on $E|U'-C'$.  
If we choose $\rho:X-C' \to [0,1]$ to be 1 near $C'$ and 0 
outside $U'$, then $J_{\rho \psi}$ extends over $C'$ as
$J$, providing the required $\omega_E$-tame complex
structure on $E$.

We finish the proof of (b) by showing that $J_0|Q_x$ is
$\omega_E$-tame for all $x \in X-C'$.  
If $g$ denotes the $J_0$-invariant, symmetric bilinear form on $Q_x$
given by $g(v,w) = \frac{1}{2}(\omega_E(v, J_0 w) +
\omega_E(w,J_0v))$, then $g(q, q) = \omega_E(q, J_0q)$, so
it suffices to show $g$ is positive definite.  
For a fixed 
$J_0$-invariant inner product on $Q_x$, let $Q_- \subset
Q_x$ be the span of all nonpositive eigenvectors of $g$. 
By $J_0$-invariance of $g$ and the background 
inner product, $Q_-$ is a $J_0$-complex subspace of $Q_x$.  
Let $\psi:Q_x \to P_x$ be a (real) 
linear transformation whose graph is a $J_\alpha$-complex
subspace of $E_x = Q_x \oplus P_x$ for some $\alpha$.  
Then $J_\psi|\graph\psi = J_\alpha|\graph\psi$
is $\omega_E$-tame.  
Thus \eqref{eq:4} is positive for all
nonzero $(q,p)$ with $p=\psi(q)$.  
The latter condition
cancels two terms, so we obtain $\omega_E(q, J_0 q) +
\omega_E(p, \psi(J_0 q))>0$. 
Since the first term is 
$g(q,q)$,  it is nonpositive on $Q_-$, implying that
$p=\psi(q)$ cannot vanish on $Q_-$ unless $q$ does. 
Hence, $\psi|Q_-$ is a monomorphism and $\dim_\complex Q_- 
\leq \dim_\complex P_x = 1$.  
It now suffices to rule out the case $\dim_\complex Q_- =1$.  
In this case, note that $J_\alpha$ and $\omega_E$ induce the same
orientation on $E_x$ and on $P_x$, hence, on $Q_x \cong
E_x/P_x$ (by $\omega_E$-orthogonality of the splitting
$E_x=Q_x \oplus P_x$). 
If $\dim_\complex Q_x = 1$, this
implies $J_0|Q_x$ is $\omega_E$-tame, so it suffices to
assume  $\dim_\complex Q_x > 1$.  
Then the function $q \mapsto g(q,q)$ realizes both positive and nonpositive
values on $Q_x -\{0\}$, so by connectedness there is a 
nonzero $q \in Q_x$ with $g(q,q) = 0$.  
It follows that 
$\omega_E$ vanishes on the $J_0$-complex line $Q_0 
\subset Q_x$ containing $q$, so it is degenerate on the 
$\omega_E$-orthogonal sum $Q_0 \oplus P_x \subset E_x$.
But this latter subspace is $J_\alpha$-complex 
since it projects to a $J_{P_x}$-complex line in $E_x/P_x$, 
contradicting the hypothesis that $J_\alpha$ is 
$\omega_E$-tame. 

To prove (c), assume $X_0 = X - \Int C$ is locally
compact, and let $\C(\xi)$ be the space of (continuous)
sections of the bundle $\xi= \End (E|X_0)$ of 
endomorphisms $E_x \to E_x$ $ (x \in X_0)$ in the 
compact-open topology.  
A subbasis for this 
topology is given by all subsets of the form $S(B,W) =
\{\sigma \in \C(\xi)|\sigma(B) \subset W\} \subset \C(\xi)$,
for $B \subset X_0$ compact and $W \subset \xi$ open.  
(See, e.g., \cite{Mu} for the space of {\em all\/} maps
$X_0 \to \xi$, then restrict to the subspace of sections $\C(\xi)$.) 
We need three basic facts about this topology: 
(1)~The evaluation map $e:X_0 \times \C(\xi) \to \xi$ (by $e(x, \sigma) = 
\sigma(x))$ is continuous.  
(This follows easily from local compactness.)  
(2)~If $S$ is any topological space and 
$\pi_{X_0} : X_0 \times S \to X_0$ is projection, then any 
section of the bundle $\pi_{X_0}^* \xi = \xi \times S
\to X_0 \times S$ has the form $(\tau, \pi_S)$ where 
$\tau:X_0 \times S \to \xi$ can be interpreted as a map $S \to \C(\xi)$.  
This correspondence gives a bijection between
continuous sections of $\pi_{X_0}^* \xi$ and continuous
maps $S \to \C(\xi)$ (cf. \cite{Mu}).  
(3)~If $X_0$ is compact, then the compact-open topology on $\C(\xi)$ equals
the $C^0$-topology induced by any fiber metric on $\xi$.
(This is clear once we observe that subtracting a section
of $\xi$ induces a fiberwise-isometric homeomorphism of 
$\xi$, so it suffices to compare neighborhoods of the $0$-section.)  

In either case (a) or (b), let $\J \subset 
\C(\xi)$ be the subspace consisting of all complex 
structures on $E|X_0$ satisfying the conclusion of the 
lemma --- that is, $(\omega_F, T)$-compatible extensions
$J$ of $J_C|\partial C$ with $T_* J = J_D$ on $D 
\cap  X_0$, and with $J$ 
$\omega_E$-tame in case (b) or (for the addendum) 
$(\omega_F,T)$-extendible along the given sequences.
For $X_0$ compact, the first step in showing the 
contractibility of $\J$ is to apply the lemma to the 
metrizable space $\tilde{X} = I \times X_0 \times \J$,
with the structures $\tilde{E}$, $\tilde{F}$, $\tilde{T}$, 
$\omega_{\tilde{F}}$ (and $\omega_{\tilde{E}}$ in case (b)) 
pulled back in the obvious way by projection $\pi_{X_0}: 
\tilde{X} \to X_0$.  
We set $\tilde{C} = ( \{0,1\} \times
X_0 \cup I \times \partial C) \times \J$ and  $\tilde{D} = 
( \{0,1\} \times X_0 \cup I \times ( D \cap X_0)) \times \J$.  
To define $J_{\tilde{C}}$ on a neighborhood 
$\tilde{U}$ of $\tilde{C}$, first consider the tautological
complex structure  $J_{\text{taut}} = ( e, \pi_\J): X_0 
\times \J \to \xi \times \J = \pi_{X_0}^* \xi =
\End (\tilde{E}|\{0\} \times X_0 \times \J)$.  
Let $h:[0, \frac{1}{2}] \times I \to I$ be continuous with 
$h(0,t) = t$, $h(s,1) = 1$ and $h(s,t) = 0$ for $s \geq t$.  
Extend to a map $[0,\frac{1}{2}] \times K \to K$, for the
collar $K \approx I \times \partial C$ in $X_0$, by 
$(h \circ (\id_{[0,1/2]} \times \pi_I), \pi_{\partial C})$
(so the $\partial C$ factor is carried along  trivially),
and let $H:[0,\frac{1}{2}] \times X_0 \to X_0$ be the 
resulting extension by $\pi_{X_0}$ outside of $K$.  
Compatibility of the collar neighborhood shows that
$(H \times \id_\J)^*(\tilde{E}|\{0\} \times X_0 
\times \J)$ can be identified with $\tilde{E} | 
([0,\frac{1}{2}] \times X_0 \times \J)$, and similarly for
$\tilde{F}$, $\omega_{\tilde{F}}$, $\tilde{T}$ (and 
$\omega_{\tilde{E}}$ and $J_{\tilde{D}}= \pi_{X_0}^* (J_D|D)$ 
when relevant), so $(H \times \id_\J)^* J_{\text{taut}}$ can be 
considered an $(\omega_{\tilde{F}}, \tilde{T})$-compatible
complex structure on $\tilde{E}|([0,\frac{1}{2}] \times X_0 \times \J)$.  
Now fix some $J \in \J$ and pull it back to 
 $\tilde{E}|([\frac{1}{2}, 1] \times X_0 \times \J)$ by the
same method, using the function $h'(s,t) = h(1-s,t)$.
Then the two complex structures agree over $\{\frac{1}{2}\} 
\times  [0, \frac{1}{2}) \times \partial C \times \J$ 
(where they are both pulled back from $J_C|\partial C$), so
together they define a complex structure $J_{\tilde{C}}$
on a neighborhood $\tilde{U}$ of $\tilde{C}$ as required.  
Construct the required local complex structures 
elsewhere on $\tilde{X}$ by (for example) pulling back the given
structures on $X_0$ by $\pi_{X_0}$. 
The hypotheses of the lemma are now 
satisfied on $\tilde{X}$ if $D$ equals $C$ or $X$. 
For the general case, 
extend $J_{\tilde D}$ over a neighborhood $\tilde V$ of $\tilde D$ by a 
procedure similar to the one above, pulling back each local structure 
on $\tilde E$ (including $J_{\tilde C}$) along the product lines of the 
given collar of $\partial D$. 
The lemma gives a section $\tilde{J}:
\tilde{X} = I \times X_0 \times \J \to \End(\tilde{E})$,
and hence a continuous map $\varphi:I \times \J \to \C(\xi)$
with image in $\J$. 
Since $\varphi |\{0\}\times \J = \id_\J$
and  $\varphi |\{1\}\times \J = J$ (as a constant map into $\J$), 
$\varphi: I \times \J \to \J$ is the required contraction.

If $X_0$ is only locally compact, we wish to show that every
map $S^m \to \J$ is homotopic to a constant map.  
But such a map can be interpreted as a section of the bundle 
$\pi_{X_0}^* \xi \to X_0\times S^m$ that is in fact a 
complex structure on $\pi_{X_0}^*(E|X_0) \to X_0 \times S^m$.  
The previous argument with $S^m$ in place of $\J$
provides the required nullhomotopy. 
\end{proof}


\begin{thebibliography}{ABKLR}

\bibitem[ABKLR]{ABKLR}
  B. Aebischer, M. Borer, M. K\"alin, Ch. Leuenberger and H. Reimann, 
  {\em Symplectic Geometry}, Progress in Math. {\bf 124}, 
  Birkh\"auser, 1994. 
  
\bibitem[A1]{A1}
  D. Auroux,
  {\em Symplectic maps to projective spaces and symplectic invariants},
  Proceedings,  G\"okova Geometry-Topology Conference 2000,
  Turk. J. Math. {\bf25} (2001), 1--42. 
  
\bibitem[A2]{A2}
  D. Auroux, 
  {\em Estimated transversality in symplectic geometry and projective maps}, 
  Proc. International KIAS Conference (Seoul 2000), 
  to appear. 
  
\bibitem[AK]{AK}
  D. Auroux and L. Katzarkov, 
  {\em The degree doubling formula for braid monodromies and
    Lefschetz pencils},
  J. Sympl. Geom., to appear.
  
\bibitem[D]{D}
  S. Donaldson,
  {\em Lefschetz fibrations in symplectic geometry}.
  Doc. Math. J. DMV, Extra Volume ICMII (1998), 309--314.
  
\bibitem[G1]{G1}
  R. Gompf,
  {\em A new construction of symplectic manifolds}.
  Ann. of Math. {\bf142} (1995), 527--595.
  
\bibitem[G2]{G2}
  R. Gompf,
  {\em The topology of symplectic manifolds}, 
  Proceedings,  G\"okova Geometry-Topology Conference 2000,
  Turk. J. Math. {\bf25} (2001), 43--59. 
  
\bibitem[G3]{G3}
  R. Gompf,
{\em Symplectic structures from Lefschetz pencils in high dimensions}, 
Geom. Topol. Monogr. {\bf 7} (2004), 267--290.
  
\bibitem[G4]{G4}
  R. Gompf,
{\em Locally holomorphic maps yield symplectic structures}, 
Comm. Anal. and Geom. {\bf 13} (2005), 511--525.

\bibitem[GS]{GS}
  R. Gompf and A. Stipsicz,
  {\em 4-manifolds and Kirby calculus}.
  Grad. Studies in Math. 20,
  Amer. Math. Soc., Providence, (1999).

\bibitem[GH]{GH}
  P. Griffiths and J. Harris,
  {\em Principles of algebraic geometry}.
  Wiley, New York, (1978).

\bibitem[Gr]{Gr}
  M. Gromov, 
  {\em Partial Differential Relations}, 
  Ergebnisse der  Mathematik, Series 3, Vol.9, Springer, Berlin, 1986.
  
\bibitem[K]{K}
  D. Kotschick,
  {\em The Seiberg-Witten invariants of symplectic four-manifolds (after
    C.H.~Taubes)},
  Seminare Bourbaki 48\`eme ann\'ee (1995-6), no.~812, Ast\'erisque {\bf241}
  (1997), 195--220.
  
\bibitem[L]{L}
  K. Lamotke,
  {\em The topology of complex projective varieties after S. Lefschetz}, 
Topology {\bf 20} (1981), 15--51. 

\bibitem[McS]{McS}
  D. McDuff and D. Salamon,
  {\em Introduction to Symplectic Topology}.
  Oxford University Press (1995).
  
\bibitem[MS]{MS}
  J. Milnor and J. Stasheff, 
  {\em Characteristic Classes}, 
  Ann. of Math. Studies 76, Princeton Univ. Press, 1974. 

\bibitem[M]{M}
  J. Moser,
  {\em On  the volume elements on a manifold}.
  Trans. Amer. Math. Soc. {\bf120} (1965), 286--294.
  
\bibitem[Mu]{Mu}
  J. Munkres, 
  {\em Topology, Second Edition}, 
  Prentice Hall, 2000. 
  
\bibitem[P]{P}
  R. Palais, 
  {\em Homotopy theory of infinite dimensional manifolds}, 
  Topology {\bf 5} (1966), 1--16.
  
\bibitem[Ta]{Ta}
  C. Taubes,
  {\em The Seiberg-Witten invariants and symplectic forms}.
  Math. Res. Letters {\bf 1} (1994), 809--822.
  
\bibitem[T]{T}
  W. Thurston,
  {\em Some simple examples of symplectic manifolds}.
  Proc. Amer. Math. Soc. {\bf55} (1976), 467--468.
  
\end{thebibliography}
\end{document}